\documentclass[twoside,11pt,reqno]{amsart}
\usepackage{amsmath,amssymb,amscd,mathrsfs}

\makeatletter

\newtheorem{thm}{Theorem}
\newtheorem{lem}[thm]{Lemma}
\newtheorem{cor}[thm]{Corollary}

\newtheorem{example}[thm]{Example}
\newtheorem{remark}[thm]{Remark}

\@addtoreset{equation}{section}

\def\row{{\operatorname{ro}}}

\def\col{{\operatorname{co}}}

\def\Col{\operatorname{Col}}
\def\Row{\operatorname{Row}}
\def\Dom{\operatorname{Dom}}
\def\Std{\operatorname{Std}}

\def\id{\operatorname{id}}

\def\bc{\bigcup}
\def\P{\Theta}
\def\D{\mathcal{D}}
\def\R{\mathcal{R}}

\def\H{\mathcal{H}}
\def\M{\mathcal{M}}
\def\T{\mathcal{T}}
\def\U{\mathcal{U}}
\def\UU{\mathscr{U}}
\def\V{\mathcal{V}}
\def\VV{\mathscr{V}}
\def\O{\mathcal{O}}

\def\C{{\mathbb C}}
\def\A{{\mathbb A}}
\def\Q{{\mathbb Q}}
\def\Z{{\mathbb Z}}
\def\N{{\mathbb N}}

\def\hom{{\operatorname{Hom}}}

\def\im{{\operatorname{Im\,}}}
\def\ker{{\operatorname{Ker\,}}}

\def\bs{\backslash}
\def\eps{{\varepsilon}}
\def\phi{{\varphi}}
\def\emptyset{{\varnothing}}

\def\bw{{{\textstyle\bigwedge}}}

\def\underbar{\mathpalette\@underbar}
\def\@underbar#1#2{\settowidth{\@tempdimb}{$#1#2$}\@tempdimb=0.8\@tempdimb
                   \ooalign{$#1#2$\crcr%
                         \hfil\rule[-.5mm]{\@tempdimb}{.4pt}\hfil}}

\newdimen\hoogte    \hoogte=11pt    % hoogte  van hokje
\newdimen\breedte   \breedte=13pt   % breedte van hokje
\newdimen\dikte     \dikte=0.5pt    % dikte lijn

\newenvironment{Young}{\begingroup
       \def\vr{\vrule height0.8\hoogte width\dikte depth 0.2\hoogte}
       \def\fbox##1{\vbox{\offinterlineskip
                    \hrule height\dikte
                    \hbox to \breedte{\vr\hfill##1\hfill\vr}
                    \hrule height\dikte}}
       \vbox\bgroup \offinterlineskip \tabskip=-\dikte \lineskip=-\dikte
            \halign\bgroup &\fbox{##\unskip}\unskip  \crcr }
       {\egroup\egroup\endgroup}
\def\diagram#1{\relax\ifmmode\vcenter{\,\begin{Young} #1\end{Young}\,}\else%
              $\vcenter{\,\begin{Young}#1\end{Young}\,}$\fi}

\begin{document}
\begin{abstract}
We derive a formula for the entries of the (unitriangular) transition matrices between 
the standard monomial and
dual canonical
bases of the irreducible 
polynomial representations of $U_q(\mathfrak{gl}_n)$ 
in terms of Kazhdan-Lusztig polynomials.
\end{abstract}
\title[Dual canonical bases]{\boldmath Dual canonical bases and Kazhdan-Lusztig polynomials}
\author{Jonathan Brundan}
\thanks{Address: Department of Mathematics, University of Oregon, Eugene, OR 97403.}
\thanks{E-mail: {\tt brundan@math.uoregon.edu}}
\thanks{Research partially supported by the NSF (grant no. DMS-0139019).}
\thanks{
{\em 2000 Subject Classification}: 17B37.}
\maketitle
\vspace{-5mm}

\iffalse
\begin{center}
{\em Dedicated to Professor Gordon James
  on the occasion of his sixtieth birthday.}
\end{center}
\fi

\section{Introduction}

In the last few years, there has been much interest in dual canonical bases
associated to quantized enveloping algebras
motivated by applications to representation theory:
in many situations 
the basis of simple modules for the Grothendieck groups of various
natural categories of modules in type $A$ can be identified with the
specialization at $q=1$ of an appropriate dual canonical basis.
For example, in \cite{BK}, we found just such an interpretation for
dual canonical bases of the irreducible polynomial representations of
$U_q(\mathfrak{gl}_n)$.
This provided the incentive to revisit the extensive literature 
about these very special modules and their bases.

The main result of the article gives an explicit formula for the
entries of the transition matrices between various 
standard monomial bases and the dual canonical basis of the irreducible
polynomial representation parametrized by a partition $\lambda
= (\lambda_1, \dots, \lambda_n)$ of $d$, 
in terms of the Kazhdan-Lusztig polynomials $P_{x,y}(t)$
associated to the symmetric group 
$S_d$. Using notation introduced later in the article,
the polynomials arising as the entries
of these matrices are of the form
$$
(-q)^{\ell(y)-\ell(x)} \sum_{z \in \D_\nu \cap S_\nu x S_\mu}(-1)^{\ell(z)+\ell(y)}
P_{z w_d, y w_d}(q^{2})
$$
for particular 
$x,y \in S_d$; see Theorem~\ref{main1} and Remark~\ref{dol}.
It is these polynomials which when evaluated at $q=1$ compute 
composition multiplicities
of 
the standard modules for the finite $W$-algebras/shifted Yangians studied in \cite{BK}. We also show that all the coefficients of these polynomials are 
non-negative integers, by relating them to the
dual canonical basis
of the quantized coordinate algebra of the group of
upper unitriangular matrices then appealing to results of Lusztig
in that setting.

The basic strategy is as follows.
Let $\mu = (\mu_1,\dots,\mu_l)$ be a composition having
transpose partition equal to $\lambda$.
Let $\V_n$ be the natural representation of 
$U_q(\mathfrak{gl}_n)$, over the field $\Q(q)$ where $q$ is an indeterminate.
By the Littlewood-Richardson rule, the space of 
$U_q(\mathfrak{gl}_n)$-module homomorphisms
$$
\xi_\mu:
\bw^{\mu_1}(\V_n)\otimes\cdots\otimes \bw^{\mu_l}(\V_n)\rightarrow
{S}^{\lambda_n}(\V_n)\otimes\cdots\otimes {S}^{\lambda_1}(\V_n)
$$
is one dimensional, and the image of any non-zero such homomorphism $\xi_\mu$ 
is the irreducible $U_q(\mathfrak{gl_n})$-module $P^\lambda(\V_n)$
of highest weight $\lambda$.
Now, the exterior and symmetric powers of $\V_n$ equipped with their
natural monomial bases
are based modules in the sense
of \cite[ch.27]{Lubook}, so by dualizing 
Lusztig's construction of tensor product
of based modules we obtain dual canonical bases for the above 
tensor products of exterior and symmetric powers. 
These bases have the remarkable property that the 
homomorphism $\xi_\mu$ 
(suitably normalized) maps dual canonical basis elements either
to dual canonical basis elements or to zero. 
In this way, we obtain the dual canonical basis of $P^\lambda(\V_n)$
($=$ the upper global crystal base of Kashiwara)
as the set of non-zero images of dual canonical basis elements of
the tensor product of exterior powers under the map $\xi_\mu$. 
Using this description, we are then able to
relate dual canonical bases directly to Kazhdan-Lusztig polynomials
using Schur-Weyl duality, following the algebraic approach 
initiated by Frenkel, Khovanov and Kirillov in \cite{FKK}.

In the main body of the article, 
we have also explained for completeness the dual argument, involving the
homomorphism
$$
\xi_\mu^*:
{S}^{\lambda_1}(\V_n)\otimes\cdots\otimes {S}^{\lambda_n}(\V_n)
\rightarrow
\bw^{\mu_l}(\V_n)\otimes\cdots\otimes \bw^{\mu_1}(\V_n)
$$
that is dual to the above map $\xi_\mu$ under certain natural pairings.
The cokernel of $\xi_\mu^*$ gives another much-studied realization of the
irreducible module $P^\lambda(\V_n)$. Again, it is the case that $\xi_\mu^*$
maps canonical basis elements either to canonical basis elements or to zero,
which makes
this point of view well-suited to relating the canonical basis
of $P^\lambda(\V_n)$ ($=$ the lower global crystal base)
to the semi-standard basis of Dipper and James \cite{DJ2}. 
In particular, we recover
the explicit formula for the transition matrix
between these bases in terms of Kazhdan-Lusztig polynomials
obtained originally by Du \cite{Du2,Du3} by a different method 
(involving the combinatorics of cells
in the symmetric group).
Along the way, we have included proofs of 
a number of related results about canonical and dual canonical bases
which are known to experts but hard to find in the literature. In
particular, in $\S$\ref{sqca}, we discuss
in some detail the dual canonical
basis of the quantized coordinate algebra of $m \times n$ matrices,
in the spirit of the work of Berenstein and Zelevinsky \cite{BZ}.
This dual canonical basis also has a 
natural representation theoretic interpretation
which does not seem to be widely known, in terms
of certain blocks 
of the categories of Harish-Chandra bimodules associated to the Lie algebras
$\mathfrak{gl}_d(\C)$.

\vspace{2mm}

\noindent
{\em Acknowledgements.}
It is a pleasure 
to thank Arkady Berenstein for numerous instructive conversations
about dual canonical bases.

\section{Combinatorics}\label{scombinatorics}

In this preliminary section, we gather together (almost) all of the 
combinatorial definitions needed later on.
Let $S_d$ denote the symmetric group acting on the
left on the set $\{1,\dots,d\}$,
with {\em basic transpositions}
$s_1,\dots,s_{d-1}$, {\em length function} $\ell$ and {\em longest element} 
$w_d$.
The following notation is quite standard:
\begin{itemize}
\item[-] $X_n$ denotes the integral weight lattice associated to the Lie algebra
$\mathfrak{gl}_n$, that is, the abelian group $\Z^n$ with standard basis $\eps_1,\dots,\eps_n$
and inner product $(.,.)$ defined by $(\eps_i,\eps_j) = \delta_{i,j}$;
\item[-] a choice of {\em simple roots} is 
given by $\eps_1-\eps_2,\dots,\eps_{n-1}-\eps_n$;
\item[-] $\geq$ is the corresponding {\em dominance ordering} on $X_n$ defined by
$\lambda \geq \mu$ if $(\lambda-\mu)$ is a sum of simple roots;
\item[-] $\Lambda_n$ and $\Lambda_n^+$ 
denote the subsets of $X_n$ consisting of all
$\lambda = (\lambda_1,\dots,\lambda_n)$
with $\lambda_1,\dots,\lambda_n \geq 0$ and with
$\lambda_1 \geq \cdots \geq \lambda_n \geq 0$, respectively;
\item[-] for a weight $\lambda\in \Lambda_n$ with
$|\lambda| := \lambda_1+\cdots+\lambda_n = d$, 
$S_\lambda$ denotes the parabolic subgroup
$S_{\lambda_1} \times \cdots \times S_{\lambda_n}$ of $S_d$
with longest element $w_\lambda$;
\item[-]
$\D_\lambda$ is the set of all {\em minimal length}
$S_\lambda \bs S_d$-coset representatives.
\end{itemize}
Letting $I_n = \{1,\dots,n\}$, 
$S_d$ also acts naturally on the right on
the set of all {\em multi-indexes}
$\alpha = (\alpha_1,\dots,\alpha_d) \in I_n^d$, so
that $(\alpha \cdot x)_i = \alpha_{xi}$ for $\alpha \in I_n^d$ and $x \in S_d$.
We write $\alpha\sim\beta$
if two multi-indexes $\alpha,\beta \in I_n^d$
lie in the same $S_d$-orbit. This is the case 
if and only if $\theta(\alpha) = \theta(\beta)$,
where $\theta(\alpha) \in \Lambda_n$ denotes 
the {\em weight} 
of $\alpha \in I_n^d$ defined from
$\theta(\alpha) = \sum_{i=1}^d \eps_{\alpha_i}$.
For $\lambda \in \Lambda_n$, let
$I_\lambda$ denote the set of all multi-indexes of weight 
$\lambda$.
There is a bijection $d:I_\lambda \rightarrow \D_\lambda$
defined
for $\alpha \in I_\lambda$ by letting $d(\alpha)$
be the unique element
of $\D_\lambda$ such that $\alpha 
\cdot d(\alpha)^{-1}$ is a weakly increasing sequence.

Assume now
that we are given weights
$\mu \in \Lambda_m$ and $\nu \in \Lambda_n$
with $|\mu|=|\nu|=d$.
The symmetric group $S_d$ acts diagonally on the right on
$I_\mu \times I_\nu$, and we let $(I_\mu \times I_\nu)/ S_d$ 
denote the set of orbits.
This set arises naturally in many different guises.
Let us recall some of the most popular.
The first involves $m \times n$ matrices 
$M = 
(m_{i,j})_{1 \leq i \leq m, 1 \leq j \leq n}$ with non-negative integer
entries. Define the {\em row} and {\em column sums} of $M$ 
to be the weights
$\row(M) = (\mu_1,\dots,\mu_m) \in \Lambda_m$ and
$\col(M) = (\nu_1,\dots,\nu_n) \in \Lambda_n$ defined from
$$
\mu_i = \sum_{j=1}^n m_{i,j}
\qquad\text{and}\qquad
\nu_j = \sum_{i=1}^m m_{i,j}.
$$
Let $\P_{\mu,\nu}$ denote the set of all such matrices
$M$ with $\row(M) = \mu$ and $\col(M) =\nu$.
Given any pair $(\alpha,\beta) \in I_\mu \times I_\nu$, we
obtain a matrix $M\in \P_{\mu,\nu}$
by letting $$
m_{i,j} = \#\{k=1,\dots,d\:|\:\alpha_k = i, \beta_k = j\}.
$$
This induces a bijective
correspondence between the sets
$(I_\mu \times I_\nu) / S_d$ 
and $\P_{\mu,\nu}$.

The second way is in terms of row standard tableaux of row shape
$\mu$ and weight $\nu$.
To introduce these, we need the notion of the
{\em row diagram} of a weight $\mu \in \Lambda_m$.
This is the diagram drawn in the positive quadrant of the $x$-$y$ plane
consisting of $\mu_1$ boxes in the first (bottom) row, \dots,
$\mu_m$ boxes in the $m$th row.
For instance, if $\mu = (5,3,4)$ its row diagram is
$$
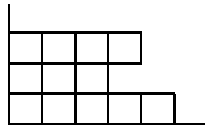
\begin{picture}(79, 38)%
\put(0,-3){\line(1,0){74}}
\put(0,-3){\line(0,1){45}}
\put(0,-3.1){\makebox(62.65,34.6){\diagram{&&&\cr&&\cr&&&&\cr}}}
\end{picture}
$$
A {\em tableau} of {row shape} $\mu$
and {weight} $\nu$
means a filling 
of the boxes of the row diagram
of $\mu$ with integers,
exactly $\nu_1$ of which are
equal to $1$,
$\nu_2$ are equal to 2,
\dots, $\nu_n$ are equal to $n$.
We sometimes use the notation $\sigma(A)$ for the {\em row shape}
$\mu$ and $\theta(A)$ for the {\em weight} $\nu$ of the
tableau $A$.
Define an equivalence relation $\sim_{\row}$ on the set of
all such tableaux 
by declaring that $A \sim_{\row} B$ if $B$ can be obtained from $A$
by permuting entries within rows.
We say that $A$ 
is {\em row standard}
if its entries are {\em weakly increasing}
along rows from left to right.
Obviously, the row standard tableaux give a set of representatives for the
$\sim_{\row}$-equivalence classes.
For instance 
\begin{equation}\label{eg}
A = \diagram{1&2&3&4\cr $2$&2&$3$\cr 1&1&2&4&$4$\cr}
\end{equation}
is a row standard tableau of row shape $(5,3,4)$ and
weight $(3,4,2,3)$.
Let $\Row(\mu,\nu)$
denote the set of all
row standard tableaux of row shape $\mu$
and weight $\nu$.
Given a tableau $A \in \Row(\mu,\nu)$, 
we obtain a matrix $M \in \P_{\mu,\nu}$ by defining
$m_{i,j}$ to be the number of entries in the $i$th row of $A$ that are equal 
to $j$.
This 
defines a bijection $\Row(\mu,\nu) \rightarrow \P_{\mu,\nu}$,
hence composing with the bijection in the previous paragraph we also
obtain a bijection between $\Row(\mu,\nu)$ and the set $(I_\mu 
\times I_\nu) / S_d$. For example, with $A$ as in (\ref{eg}),
the corresponding matrix $M \in \P_{\mu,\nu}$ is the matrix
$$
\left(
\begin{array}{cccc}
2&1&0&2\\
0&2&1&0\\
1&1&1&1
\end{array}
\right)
$$
and a representative $(\alpha,\beta) \in I_\mu \times I_\nu$
for the corresponding orbit is given by setting
$\alpha = (3,3,3,3,2,2,2,1,1,1,1,1)$ and $\beta = (1,2,3,4,2,2,3,1,1,2,4,4)$.

The third way 
involves the set $\D^+_{\nu,\mu}$ of {\em maximal length}
distinguished $(S_\nu, S_\mu)$-double coset representatives in the
symmetric group $S_d$. 
We just explain how to define a bijection between
$\Row(\mu,\nu)$ and $\D_{\nu,\mu}^+$.
Given any tableau $A$ of row 
shape $\mu$
and weight $\nu$, define a sequence
$\rho(A) \in I_\nu$ by {\em row reading} the entries of $A$ along rows from 
left to right starting from the top row;
for example, if $A$ is as in
(\ref{eg}) then $\rho(A)$ 
is the multi-index $\beta$ from the end of the previous paragraph.
Recalling the bijection $d:I_\nu \rightarrow D_\nu$ from the opening paragraph,
the map $A \mapsto d(\rho(A)) w_d$
defines a bijection between the set $\Row(\mu,\nu)$
and the set $D_{\nu,\mu}^+$.
Moreover,
$D_\nu \cap S_\nu d(\rho(A))w_d  S_\mu
=
\{w_\nu d(\rho(B))w_d\:|\:B \sim_\row A\}$.
For a proof of a similar statement, 
see \cite[1.7]{DJ} or \cite[4.4]{Mat}.

There is a fourth way which is much more subtle than the ones discussed so 
far involving the Robinson-Schensted-Knuth correspondence; 
see \cite[$\S$4.1]{fulton}.
In this article, we actually only need a very special case of this fundamental 
bijection. To explain it, we must first
introduce the 
notion of a {column strict tableau}. 
Suppose now
that $\mu \in \Lambda_l$, $\nu \in \Lambda_n$ satisfy $|\mu|=|\nu|=d$.
The mirror image of the row diagram of $\mu$
in the line $y=x$ gives the {\em column diagram} of $\mu$.
Thus, the column diagram has $\mu_1$ boxes in the 
first (leftmost) column, \dots, $\mu_l$ boxes in the $l$th
column.
A {tableau} of {\em column shape} $\mu \in \Lambda_l$ 
and weight $\nu \in \Lambda_n$ means a filling of the boxes
of the {column diagram} 
of $\mu$ with integers, 
exactly $\nu_j$ of which are equal to $j$ for each $j=1,\dots,n$.
Call such a tableau
{\em column strict} if its entries are
{\em strictly increasing} along columns from bottom to top.
Let $\Col(\mu,\nu)$ denote the
set of all column strict tableaux of column shape
$\mu$ and weight $\nu$.
Observe that the mirror image
in the line $y=x$ of a tableau $A$ of column shape $\mu$ defines a tableau
$A'$ of row shape $\mu$. 
This is a useful trick for carrying over the earlier definitions  to
the present setting. For instance, we 
write $A \sim_{\col} B$ if $A' \sim_{\row} B'$.
The next definition breaks the symmetry:
define the {\em column reading} $\gamma(A)$ to be the multi-index
obtained by reading the entries of $A$ 
along columns from top to bottom starting from the leftmost column.
This is related to the row reading of $A'$ by the equation
$\gamma(A) = \rho(A') \cdot w_d$.

Assuming all parts of  the composition $\mu$ are $\leq m$, 
let
$\lambda = \mu' \in \Lambda_m^+$ be the {\em conjugate partition},
so $\lambda_i$ 
is the number of boxes in the $i$th row of the column
diagram of $\mu$.
Let $\Dom(\lambda,\nu)$ denote the familiar set
of all {\em standard tableaux} of row shape $\lambda$ and
weight $\nu$, that is, the tableaux in $\Row(\lambda,\nu)$ that are
also column strict. 
(The unfamiliar symbol $\Dom$ here stands for ``dominant'' following
the language used in \cite{BK}.)
For a multi-index $\alpha \in I_n^d$,
let $P(\alpha)$
denote the image of the word
$\alpha_1\alpha_2\cdots \alpha_d$ under the Robinson-Schensted correspondence;
see e.g. \cite[$\S$4.1]{fulton}. Thus, 
$P(\alpha)$ is the standard tableau 
$\varnothing \leftarrow \alpha_1 \leftarrow\cdots \leftarrow \alpha_d$, 
where $\leftarrow$ denotes row insertion as in \cite[$\S$1.1]{fulton}.
Still writing $\lambda = \mu'$, define
\begin{equation}
\Std(\mu,\nu) = \{A \in \Col(\mu,\nu)\:|\:P(\gamma(A))
\text{ is of row shape }\lambda\}.
\end{equation}
We refer to elements of $\Std(\mu,\nu)$ as
{\em standard tableaux} of column shape $\mu$ and weight $\nu$.
In the special case
$\mu$ is itself a partition,
it is easy to see from the definition of the Robinson-Schensted
map that 
$\Std(\mu,\nu)$ is the set of all tableaux in $\Col(\mu,\nu)$
that are also row standard, i.e. $\Std(\mu,\nu) = \Dom(\lambda,\nu)$.
So the double meaning of the phrase ``standard tableaux'' is unambiguous.
In general, by a result of Lascoux and Sch\"utzenberger 
\cite{LS2}, the {\em rectification map}
\begin{equation}\label{rectify}
R:\Std(\mu,\nu) \rightarrow \Dom(\lambda,\nu), \qquad
A \mapsto P(\gamma(A))
\end{equation}
is a bijection; see also \cite[$\S$A.5]{fulton}.
In the special case that $\mu$ is a partition, the map $R$
is just the identity map. In general, 
$R$ can be computed by repeatedly using {\em jeu de taquin}
to permute adjacent columns of different lengths;
see \cite[$\S$4]{LT} for an example.

In proofs, we will use a rather different characterization of the set $\Std(\mu,\nu)$
and the rectification map in terms of crystals.
To recall this, define a crystal 
$(I_n^d, \tilde e_i, \tilde f_i, \eps_i, \varphi_i, \theta)$
in the sense of Kashiwara \cite{Ka}
with underlying set
$I_n^d$ as follows.
For $i=1,\dots,n-1$, 
define the {\em $i$-signature} 
$(\sigma_1,\dots,\sigma_d)$
of $\alpha \in I_n^d$
by\begin{equation*}
\sigma_j = \left\{
\begin{array}{ll}
+&\hbox{if $\alpha_j = i$,}\\
-&\hbox{if $\alpha_j = i+1$,}\\
0&\hbox{otherwise.}
\end{array}\right.
\end{equation*}
From this the {\em reduced $i$-signature} is computed
by successively replacing 
subsequences of the form $-+$ (possibly separated by $0$'s)
in the signature with $0$'s
until no $-$ appears to the left of 
a $+$.
Let $\delta_j$ denote the 
$d$-tuple $(0,\dots,0,1,0,\dots,0)$
where $1$ appears in the $j$th place.
Now define
\begin{align*}
\tilde e_i(\alpha) &:= \left\{
\begin{array}{ll}
\emptyset&\hbox{if there are no $-$'s in the reduced $i$-signature},\\
\alpha - \delta_j&\hbox{if the leftmost $-$ is in position $j$;}
\end{array}\right.\\
\tilde f_i(\alpha) &:= \left\{
\begin{array}{ll}
\emptyset&\hbox{if there are no $+$'s in the reduced $i$-signature},\\
\alpha + \delta_j&\hbox{if the rightmost $+$ is in position $j$;}
\end{array}\right.\\
\eps_i(\alpha) &= \hbox{the total number of $-$'s in the
reduced $i$-signature},\\
\varphi_i(\alpha) &= \hbox{the total number of $+$'s in the reduced 
$i$-signature}.
\end{align*}
Recalling that $\theta(\alpha)$ denotes the weight
of $\alpha \in I_n^d$,
this completes the definition of the crystal
$(I_n^d, \tilde e_i, \tilde f_i, \eps_i, \varphi_i, \theta)$.
It is just the $d$-fold tensor product of the usual crystal associated
to the natural $\mathfrak{gl}_n$-module, except that we have parametrized it 
from right to left rather than from left to right.

In this paragraph, we write $\bigcup$ as shorthand for the union over all
$\nu \in \Lambda_n$, and assume in addition that $m \leq n$.
The row reading $\rho$ resp. the column reading $\gamma$
identifies the set $\bigcup \Row(\lambda,\nu)$
resp. $\bigcup \Col(\mu,\nu)$
 with a subcrystal of
$I_n^d$. This defines new crystals
$(\bigcup \Row(\lambda,\nu), \tilde e_i, \tilde f_i, \eps_i, \phi_i,
\theta)$ and
$(\bigcup \Col(\mu,
\nu), \tilde e_i, \tilde f_i, \eps_i, \phi_i,
\theta)$.
It is well known that the map
$A \mapsto P(\gamma(A))$ arising from the Robinson-Schensted correspondence
commutes in the strict sense with the crystal operators
$\tilde e_i, \tilde f_i$.
Moreover, $\bigcup \Dom(\lambda,\nu)$ is a subcrystal of
$\bigcup \Row(\lambda,\nu)$, indeed, it is precisely
the connected component of
$\bigcup\Row(\lambda,\nu)$ generated by the unique tableau 
$B \in \Dom(\lambda,\lambda)$, i.e. 
the tableau with all entries in its $i$th
row equal to $i$. Since $R$ necessarily 
maps the unique element $A \in \Std(\mu,\lambda)$ to this tableau 
$B$, we deduce that $\bigcup\Std(\mu,\nu)$ is the connected component
of $\bigcup \Col(\mu,\nu)$ generated by $A$,
and the rectification map $R:\bigcup \Std(\mu,\nu) \rightarrow
\bigcup\Dom(\lambda,\nu)$ is an isomorphism of crystals.
In this way, we obtain various different realizations
of the usual highest weight
crystal associated to the partition $\lambda$,
one for each composition $\mu$
with $\mu' = \lambda$.
The standard realization from \cite{KN} is the one when
$\mu$ is itself a partition.

Finally, 
we say a few words about the Bruhat ordering.
Let $\leq$ denote the {\em opposite} of the usual Bruhat ordering on $S_d$,
e.g. $w_d \leq 1$. 
This restricts to a partial
ordering on the subset $\D_{\nu,\mu}^+$, 
for $\mu \in \Lambda_m, \nu\in\Lambda_n$ with
$|\mu|=|\nu|=d$ as before. Hence
using the above bijections,
we get partial orderings also denoted $\leq$
on each of sets $(I_\mu \times I_\nu) / S_d, 
\P_{\mu,\nu}$ and $\Row(\mu,\nu)$.
We want to record several
equivalent ways of defining these
partial orders directly; see \cite[1.2]{DJ} or 
\cite[3.8]{Mat}
for proofs of essentially the same statements, which are apparently
due originally to Ehresmann.
Suppose first that we are given tableaux $A$ and $B$. 
Write $A \downarrow B$ 
if there exists
an entry $x$ in the $i$th row 
and an entry
$y$ in the $j$th row of $A$
with $i < j$ and $x < y$ 
such that $B$ is obtained from 
$A$ by swapping the entries $x$ and $y$.
For example,
$$
\diagram{1&2&5\cr7&7\cr3&3&5\cr}
\quad
\downarrow
\quad
\diagram{1&2&3\cr7&7\cr3&5&5\cr}
\quad
\downarrow
\quad
\diagram{1&2&3\cr7&3\cr7&5&5\cr}
$$
Then, $A \geq B$ in the Bruhat ordering on $\Row(\mu,\nu)$ 
if and only if there exist
tableaux $C_1,\dots,C_r$  such that 
$A \sim_{\row} 
C_1 \downarrow \cdots \downarrow C_r \sim_{\row} B$.
Given $A \in \Row(\mu,\nu)$, let 
$A_{\leq i}$ denote the
tableau obtained from $A$ by deleting all boxes in 
rows higher than the $i$th row, and let
$A^{\leq j}$ denote the tableau obtained from $A$ by deleting all
boxes containing entries greater than $j$.
The following are equivalent
for $A, B \in \Row(\mu,\nu)$:
\begin{itemize}
\item[(i)] $A \leq B$ in the Bruhat ordering on $\Row(\mu,\nu)$;
\item[(ii)] $\theta(A_{\leq i}) \leq \theta(B_{\leq i})$ 
in the dominance ordering on $\Lambda_n$ for 
all $i=1,\dots,m$ (recall $\theta$ denotes weight);
\item[(iii)] $\sigma(A^{\leq j}) \leq \sigma(B^{\leq j})$ in the dominance
ordering on $\Lambda_m$
for all $j=1,\dots,n$ (recall $\sigma$ denotes row shape).
\end{itemize}
From (ii) or (iii), one easily deduces the well known 
direct description of
the Bruhat order on the set
$\P_{\mu,\nu}$ itself: for $M, N \in \P_{\mu,\nu}$,
we have that $M \leq N$ if and only if
$$
\sum_{i=1}^s \sum_{j=1}^t m_{i,j} \leq \sum_{i=1}^s \sum_{j=1}^t n_{i,j}
$$
for all $s=1,\dots,m$ and $t=1,\dots,n$.

We will also need the Bruhat ordering $\leq'$ on
the set $\Col(\mu,\nu)$. This can be defined
simply by $A \leq' B$ if $A' \geq B'$;
equivalently, $d(\gamma(A)) \geq d(\gamma(B))$.
In the special case that $\mu$ is
a partition and $\lambda = \mu'$, 
we have now defined two partial orders
$\leq'$ and $\leq$ on the set
$\Std(\mu,\nu) = \Dom(\lambda,\nu)$,
via its natural embeddings into $\Col(\mu,\nu)$ and $\Row(\lambda,\nu)$,
respectively.
The following lemma shows that these two partial orders
coincide.

\begin{lem}\label{bor2}
For $A, B \in \Dom(\lambda,\nu)$, 
we have that $A\leq B$ if and only if $A \leq' B$.
\end{lem}

\begin{proof}
Let $A, B \in \Dom(\lambda,\nu)$.
Since $A$ is standard, 
$\sigma((A')^{\leq j}) = \sigma((A^{\leq j})')
=
\sigma(A^{\leq j})'$,
and similarly for $B$. 
By the third 
equivalent definition of the Bruhat ordering on $\Row(\lambda,\nu)$ 
above, we know that
$A \leq B$ if and only if
$\sigma(A^{\leq j}) \leq \sigma(B^{\leq j})$
for all $j=1,\dots,n$.
Since conjugation is order reversing on partitions, this is equivalent to
$\sigma(A^{\leq j})' \geq \sigma(B^{\leq j})'$ for all $j=1,\dots,n$,
i.e. $\sigma((A')^{\leq j}) \geq \sigma((B')^{\leq j})$.
This is the statement that $A' \geq B'$, hence $A \leq' B$.
\end{proof}

\section{Quantized enveloping algebras}\label{srmatrices}

In this section, we recall the definition of the quantized enveloping algebra
$\U_n = U_q(\mathfrak{gl}_n)$, following
\cite{Lubook}.
We will work over the field $\mathbb{Q}(q)$ where
$q$ is an indeterminate.
An additive map $f:V \rightarrow W$ between $\Q(q)$-vectors
spaces is called {\em antilinear} 
if $f(cv) = \bar c f(v)$ for all $c \in \Q(q), v \in V$,
where $-:\Q(q)\rightarrow\Q(q)$ is the field automorphism with
$\overline{q} = q^{-1}$.
Also the quantum integer associated to $n \in \N$ is $[n] = (q^n - q^{-n})/(q-q^{-1})$ and the quantum factorial is
$[n]! = [n][n-1] \cdots [2][1]$.

By definition,
$\U_n$ is the $\Q(q)$-algebra on generators
$E_i, F_i\:(i=1,\dots,n-1)$ and $K_i, K_i^{-1}\:(i =1,\dots,n)$
subject to relations
\begin{align*}
K_iK_i^{-1} &= K_i^{-1}K_i = 1,&E_iE_j &= E_jE_i&\hbox{if $|i-j|>1$},\\
K_iK_j &= K_jK_i,&\!\!\!\!E_i^2E_j  + E_jE_i^2 &= [2]E_iE_jE_i&\hbox{if $|i-j|=1$},\\
K_iE_jK_i^{-1} &= q^{(\eps_i, \eps_j - \eps_{j+1})} E_j,&F_iF_j &= F_jF_i&\hbox{if $|i-j| > 1$,}\\
\qquad\quad K_iF_jK_i^{-1} &= q^{(\eps_i, \eps_{j+1}-\eps_j)} F_j,&\!\!\!\!F_i^2F_j+ F_jF_i^2 &=   [2]F_iF_jF_i&\hbox{if $|i-j|=1$,}
\end{align*}
\begin{equation*}
\quad 
E_iF_j - F_jE_i = \delta_{i, j} \frac{K_{i, i+1} - K_{i+1, i}}{q - q^{-1}}.
\end{equation*}
Here, $K_{i, j}$ denotes
$K_iK_j^{-1}$. 
In this article, we will always view $\U_n$ as a Hopf algebra with 
counit $\eps:\U_n \rightarrow \Q(q)$ defined by 
$\eps(E_i) = 0, \eps(F_i) = 0$ and
$\eps(K_i) = 1$,
and comultiplication
$\Delta: \U_n \rightarrow \U_n \otimes \U_n$ defined by
$$
\Delta(E_i) = 1 \otimes E_i + E_i \otimes K_{i+1,i},\quad\!
\Delta(F_i) = K_{i, i+1}\otimes F_i + F_i\otimes 1,\quad\!
\Delta(K_i) = K_i \otimes K_i.
$$
In the language of \cite{KaG}, 
this comultiplication is adapted to taking tensor products
of lower crystal bases at $q=0$ and upper crystal bases at $q=\infty$.
With only minor adjustments,
it would also be perfectly possible to use throughout the article
the comultiplication 
$\widetilde\Delta:\U_n \rightarrow \U_n \otimes \U_n$
from \cite{Lubook}.
We just note for comparison that $\widetilde\Delta$ is defined by
$\widetilde\Delta = 
(\tau\otimes \tau)\circ\Delta\circ\tau$, where $\tau:\U_n \rightarrow \U_n$
is the algebra antiautomorphism defined by $\tau(E_i) = F_i, \tau(F_i) = E_i$
and $\tau(K_i) = K_i$, and it is adapted to taking tensor products of
lower crystal bases at $q=\infty$ and upper crystal bases at $q=0$.

All $\U_n$-modules encountered in this article
will be
{\em polynomial representations}, meaning 
$\U_n$-modules $V$ satisfying 
$V = \bigoplus_{\lambda \in \Lambda_n} V_\lambda$ where 
$V_\lambda$ denotes the
{\em $\lambda$-weight space}
\begin{equation*}
V_\lambda = \{v \in V\:|\: K_i v = q^{(\lambda,\eps_i)} v\hbox{ for all $i=1,\dots,n$}\}.
\end{equation*}
Direct sums, tensor products and subquotients 
of polynomial representations are again polynomial.
Moreover, 
the category of all polynomial representations of $\U_n$
is a {\em braided tensor category}, with braiding
isomorphism $\R_{V,W}:V \otimes W \rightarrow W \otimes V$
defined like in \cite[32.1.5]{Lubook}.
To review this definition in a little more detail, let $\Theta
= \sum_{0 \leq \lambda \in X_n} \Theta_\lambda$ 
be the {\em quasi-$R$-matrix} defined as in \cite[4.1.2]{Lubook}, 
but using our comultiplication $\Delta$ instead of the comultiplication
$\widetilde\Delta$ used there.
More precisely,
$\Theta = (\tau \otimes \tau)(\widetilde\Theta^{-1})$
where $\widetilde\Theta$
is exactly Lusztig's {quasi-$R$-matrix} from \cite[4.1.2]{Lubook}.
It is an element of a certain
completion $(\U_n \otimes \U_n)^\wedge$ of the algebra 
$\U_n \otimes \U_n$, with $\Theta_0 = 1$ and
$\Theta_\lambda \in \U_\lambda^+ \otimes \U_\lambda^-$ for each $\lambda$, where
$\U_\lambda^+$ resp. $\U_\lambda^-$ denotes the $\pm\lambda$-weight space
of the positive part $\U_n^+$ resp. the negative part $\U_n^-$ of $\U_n$.
For polynomial representations $V$ and $W$,
all but finitely many $\Theta_\lambda$ act as zero
on any given vector $v \otimes w \in V \otimes W$,
by weight considerations.
Hence it makes sense to view $\Theta$ as 
an invertible operator on $V \otimes W$. 
The braiding $\R_{V,W}:V \otimes W \rightarrow W \otimes V$
can now be defined to be the map
$\R_{V,W} = \Theta\circ f \circ P$
where $f:W \otimes V \rightarrow W \otimes V$ is the map
$w \otimes v \mapsto q^{(\lambda,\mu)} w \otimes v$
for $v,w$ of weights $\lambda,\mu$, respectively,
and $P:V \otimes W\rightarrow W \otimes V$ is the permutation
operator $v \otimes w \mapsto w \otimes v$.
Suppose more generally that
$V_1,\dots,V_d$ are all polynomial representations.
For $1 \leq i < d$, let
$$
\R_i: V_1 \otimes\cdots\otimes V_i \otimes V_{i+1}\otimes\cdots\otimes V_d
\rightarrow 
V_1 \otimes\cdots\otimes V_{i+1} \otimes V_{i}\otimes\cdots\otimes V_d
$$
denote the $\U_n$-module isomorphism
$\R_{V_i,V_{i+1}}$ acting on the $i$th and 
$(i+1)$th tensor positions.
For a permutation $w \in S_d$, we obtain a well-defined map
$$
\R_w:V_1 \otimes \cdots \otimes
V_d \rightarrow V_{w^{-1} 1} \otimes\cdots \otimes V_{w^{-1} d}
$$
by setting $\R_w = \R_{i_1} \circ\cdots\circ \R_{i_\ell}$ 
if $s_{i_1} s_{i_2} \cdots s_{i_\ell}$ is a reduced expression for $w$.

In the remainder of the section, we want to discuss some properties
of bar involutions.
The bar involution on $\U_n$
is the unique antilinear automorphism such that
$\overline{E}_i =E_i$,
$\overline{F}_i = F_i$ and 
$\overline{K}_i = K_i^{-1}$.
We say that a $\U_n$-module $V$ possesses a {\em compatible bar involution}
if it is equipped with an antilinear involution $-:V \rightarrow V$
such that $\overline{u v} = \overline{u}\,\overline{v}$ for each $u \in \U_n$
and $v \in V$. Suppose $V$ and $W$ are polynomial
$\U_n$-modules with compatible
bar involutions. Following \cite[27.3.1]{Lubook}, 
there is a canonical way to define 
a compatible bar involution on the tensor product $V \otimes W$:
given $v \in V$ and $w \in W$ we set
\begin{equation}\label{bardef}
\overline{v \otimes w} = \Theta (\overline{v} \otimes \overline{w}).
\end{equation}
More generally, given polynomial $\U_n$-modules $V_1,\dots,V_d$ each
possessing a compatible bar involution, there is a compatible
bar involution on $V_1 \otimes \cdots \otimes V_d$ defined as follows:
pick any $1 \leq k < d$ then set
\begin{equation*}
\overline{v_1 \otimes\cdots\otimes v_d}
= 
\Theta ((\overline{v_1\otimes\cdots\otimes v_k}) \otimes
(\overline{v_{k+1}\otimes\cdots\otimes v_d}))
\end{equation*}
where the bar involutions on the right hand side are defined inductively.
By \cite[27.3.6]{Lubook}, this definition is independent of the 
particular choice of $k$.
Alternatively, in terms of the braiding,
the bar involution on $V_1 \otimes\cdots\otimes V_d$
satisfies
\begin{equation}\label{eq18}
\overline{v_1 \otimes \cdots \otimes v_d} = 
q^{-\sum_{i < j} (\lambda_i,\lambda_j)}
\R_{w_d}  (\overline{v_d}\otimes\cdots\otimes\overline{v_1})
\end{equation}
if $v_i$ is of weight $\lambda_i$,
recalling that $w_d$ denotes
the longest element of $S_d$.
Also,
\begin{equation}\label{compr}
\overline{\R_w (v_1 \otimes \cdots \otimes v_d)}
= \R_{w^{-1}}^{-1} 
(\overline{v_1 \otimes \cdots \otimes v_d})
\end{equation}
for any $w \in S_d$ and $v_i \in V_i$;
the proof of this reduces easily to the case $\ell(w) = 1$ which follows 
using the identity $\Theta^{-1} = \overline{\Theta}$ from
\cite[4.1.3]{Lubook}.

We say that $A$ is a {\em polynomial $\U_n$-algebra} 
if $A$ is a polynomial $\U_n$-module and an 
associative algebra, with identity element $1_A$ and
multiplication $\mu_A:A \otimes A \rightarrow A$,
such that $u 1_A = \eps(u) 1_A$ 
and 
$u (xx') = \mu_A(\Delta(u) (x\otimes x'))$ for each $u \in \U_n, x,x' \in A$.
Given two polynomial $\U_n$-algebras $A$ and $B$,
the tensor product $A \otimes B$ is a polynomial $\U_n$-module; we
make it into a polynomial $\U_n$-algebra by defining the multiplication
$\mu_{A \otimes B}:A \otimes B \otimes A \otimes B \rightarrow A \otimes B$
from
$
\mu_{A\otimes B} = 
(\mu_A \otimes \mu_B) \circ (\id_A \otimes \R_{B,A} \otimes \id_B).
$
It is well known that this multiplication is associative.
More generally, given polynomial $\U_n$-algebras $A_1,\dots,A_d$,
we make the tensor product $A_1 \otimes \cdots \otimes A_d$
into a polynomial $\U_n$-algebra by iterating this construction.
Explicitly, the multiplication is the map
$(\mu_{A_1}\otimes
\cdots \otimes \mu_{A_d})
\circ \R_{w}:A_1\otimes\cdots\otimes A_d \otimes A_1\otimes\cdots\otimes A_d
\rightarrow A_1\otimes\cdots\otimes A_d$
where 
$w:(1,2,\dots,d,d+1,d+2,\dots,2d) \mapsto
(1,3,\dots,2d-1,2,4,\dots,2d)$.

\begin{lem}\label{barprop}
Suppose that $A_1,\dots,A_d$ are polynomial $\U_n$-algebras
equipped with compatible bar involutions
such that
$\mu_{A_i} (\overline{x_i \otimes y_i})
= \overline{x_i y_i}$
for each $i$ and $x_i,y_i \in A_i$.
View the tensor product $A_1 \otimes \cdots \otimes A_d$
as a polynomial $\U_n$-algebra equipped with a compatible bar involution
by the above constructions.
Let $*$ denote the {twisted multiplication}
on $A_1 \otimes\cdots\otimes A_d$ defined by the map
$((\mu_{A_1}\circ \R_{A_1, A_1}) \otimes\cdots\otimes
(\mu_{A_d}\circ \R_{A_d,A_d}))
\circ \R_w
$
where $w:(1,2,\dots,d,d+1,d+2,\dots,2d)
\mapsto (1,3,\dots,2d-1,2,4,\dots,2d)$.
Then,
\begin{equation*}
\overline{(x_1 \otimes\cdots\otimes x_d)(y_1\otimes\cdots\otimes y_d)} 
=
q^{-(\lambda,\mu)}
(\overline{y_1\otimes\cdots\otimes y_d} )
*
(\overline{x_1 \otimes\cdots\otimes x_d\phantom{y_1\!\!\!\!\!\!}})
\end{equation*}
for $x_i,y_i \in A_i$ such that $x_1\otimes\cdots\otimes x_d$ is
of weight $\lambda$ and $y_1\otimes\cdots\otimes y_d$ is of weight $\mu$.
\end{lem}

\begin{proof}
Using (\ref{compr}) and the definitions, we have that
\begin{align*}
&\overline{(x_1 \otimes\cdots\otimes x_d)(y_1\otimes\cdots\otimes y_d)} \\
&\qquad=
\overline{
((\mu_{A_1}\otimes\cdots\otimes \mu_{A_d})
\circ
\R_w)(x_1 \otimes\cdots\otimes
 x_d \otimes y_1 \otimes\cdots\otimes y_d)}\\
&\qquad=((\mu_{A_1}\otimes\cdots\otimes \mu_{A_d})
\circ \R_{w^{-1}}^{-1})
(\overline{x_1 \otimes\cdots\otimes x_d \otimes y_1 \otimes\cdots\otimes y_d})\\
&\qquad=q^{-(\lambda,\mu)}
((\mu_{A_1}\otimes\cdots\otimes \mu_{A_d})\circ \R_{w^{-1}}^{-1}\circ \R_{v}) 
(\overline{y_1\otimes\cdots\otimes y_d} \otimes 
\overline{x_1\otimes\cdots\otimes x_d\phantom{y_1\!\!\!\!\!\!}})
\end{align*}
where $w$ is as in the statement of the lemma
and $v = (1\;d+1)(2\;d+2)\cdots (d\;2d)$.
Now the proof is completed by observing that
$wvw^{-1} = (1\;2)(3\;4)\cdots(2d-1\;2d)$ and
then checking that lengths add correctly so that
$\R_v = \R_{w^{-1}} \R_{wvw^{-1}} \R_{w}$.
\end{proof}

\section{Kazhdan-Lusztig polynomials}\label{skl}

The next job is to review the
definition of the parabolic Kazhdan-Lusztig
polynomials associated to the symmetric group $S_d$.
Let $\H_d$ denote the corresponding Hecke algebra.
By definition, this is the
$\Q(q)$-algebra with basis $\{H_x\:|\:x \in S_d\}$ and multiplication defined
by the rules that $H_{x} H_y = H_{xy}$ if $\ell(xy) = \ell(x)+\ell(y)$, and
\begin{equation}\label{qrel}
H_i^2 = 1-(q-q^{-1})H_i,
\end{equation}
where we write $H_i = H_{s_i}$ for short.
Take any weight $\lambda \in \Lambda_n$ with $|\lambda| = d$.
Corresponding to the parabolic subgroup $S_\lambda$ of $S_d$, we have the
parabolic
subalgebra $\H_\lambda$ of $\H_d$ spanned by  $\{H_x\:|\:x \in S_\lambda\}$.
Let $1_{\H_\lambda}$ denote the one dimensional right $\H_\lambda$-module
spanned by a vector $1_\lambda$ such that $1_\lambda H_i = q^{-1} 1_\lambda$
for each $H_i \in \H_\lambda$.
Form the induced module
\begin{equation}\label{permdef}
\M^\lambda = 1_{\H_\lambda} \otimes_{\H_\lambda} \H_d.
\end{equation}
This has a natural basis $\{M_x\:|\:x \in \D_\lambda\}$
defined from $M_x = 1_\lambda \otimes H_{x}$.
Now we can introduce the two 
families of parabolic Kazhdan-Lusztig polynomials, following
\cite[$\S$3]{soergel} closely.
We need the
bar involution on $\H_d$, that is, the unique
antilinear automorphism of $\H_d$
such that 
$\overline{H_w} = H_{w^{-1}}^{-1}$
for each $w \in S_d$; 
in particular, $\overline{H_i} = H_i - (q-q^{-1})$.
There is an induced bar involution on $\M^\lambda$, with
$\overline{M_x} = 1_\lambda \otimes \overline{H_x}$ for each
$x \in \D_\lambda$.
By \cite[3.1,3.5]{soergel}, there are unique bar invariant elements
$\underbar{\tilde M}_x, \underbar{M}_x 
\in \M^\lambda$ for each $x \in \D_\lambda$
such that
$$
\underbar{\tilde M}_x \in {M}_x + \sum_{y \in \D_\lambda} q^{-1}\Z[q^{-1}] {M}_y,
\qquad
\underbar{M}_x \in M_x + \sum_{y \in \D_\lambda} q\Z[q] M_y.
$$
In Soergel's notation, we have that
\begin{equation}\label{ekl1}
\underbar{\tilde M}_y = \sum_{x \in \D_\lambda} (-1)^{\ell(x)+\ell(y)} n_{x,y}(q^{-1}) 
M_x,\qquad
\underbar{M}_y = \sum_{x \in \D_\lambda} m_{x,y}(q) M_x
\end{equation}
for polynomials $n_{x,y}(q), m_{x,y}(q) \in \Z[q]$
which up to a shift are the usual parabolic Kazhdan-Lusztig polynomials
of \cite{KL,deodhar}; see \cite[3.2]{soergel} for the precise identification. Recalling that $\leq$ is the opposite of the usual Bruhat ordering on $S_d$, we have that $n_{x,x}(q) = m_{x,x}(q) = 1$ and 
$n_{x,y}(q) = m_{x,y}(q) = 0$ unless $x \geq y$.

We now want to 
review a 
completely different approach to the
construction of these polynomials
 involving
the quantized enveloping algebra
$\U_n = U_q(\mathfrak{gl}_{n})$ from $\S$\ref{srmatrices}
in place of the Hecke algebra $\H_d$.
The coincidence here is well explained algebraically by
Schur-Weyl duality, and that is the point of view we will take.
The exposition in the remainder of the section 
is equivalent to 
that of \cite{FKK}, which we believe is the first place that this 
elementary approach
appeared explicitly in the literature.
There is also an older geometric explanation 
which relies on the local isomorphism
between Schubert varieties
and the varieties arising from representations of quivers 
in type $A$ from \cite{Zel}; see \cite{GL}.
To start with, let $\V_n$ denote the 
natural $\U_n$-module, that is, the polynomial representation
on basis $\{v_i\:|\:i=1,\dots,n\}$ with
action defined by
\begin{equation*}
K_i v_j = q^{(\eps_i,\eps_j)} v_j,
\qquad
E_i v_j = \delta_{i+1,j} v_i,
\qquad
F_i v_j = \delta_{i,j} v_{i+1}.
\end{equation*}
The {tensor algebra} $T(\V_n) = \bigoplus_{d \geq 0} T^d(\V_n)$
is a polynomial $\U_n$-algebra in the sense of $\S$\ref{srmatrices}.
The $\U_n$-module $\V_n$
possesses compatible bar involution
defined simply by
$\overline{v_i} = v_i$ 
for each $i =1,\dots,n$.
By the tensor product construction from $\S$\ref{srmatrices}, we get induced
a compatible bar involution on each $T^d(\V_n)$, hence on
the tensor algebra $T(\V_n)$ itself.
The bar involution on $\V_n \otimes \V_n$ satisfies
\begin{equation}\label{barinv}
\overline{v_i \otimes v_j} =
\left\{
\begin{array}{ll}
v_i \otimes v_j& \hbox{if $i \leq j$,}\\
v_i \otimes v_j + (q-q^{-1}) v_j \otimes v_i&\hbox{if $i > j$.}
\end{array}\right.
\end{equation}
This can be seen as follows:
if $i \leq j$ all $\Theta_\lambda$
except for $\Theta_0$ annihilate $v_i \otimes v_j$ by weight considerations
hence $\overline{v_i \otimes v_j} = v_i \otimes v_j$ in these cases;
then for $i > j$ one applies $F_{i-1}F_{i-2}\cdots F_j$
to both sides of the identity 
$\overline{v_j \otimes v_j} = v_j \otimes v_j$
to deduce the formula in these cases too.

Combining (\ref{barinv}) with (\ref{eq18}), one checks that the
inverse braiding $\R^{-1}_{\V_n, \V_n}$ satisfies the quadratic relation
(\ref{qrel}).
Hence, there is a well-defined
{\em right action} of the Hecke algebra $\H_d$
on $T^d(\V_n)$
defined from $v H_w = \R_{w}^{-1} (v)$
for $v \in T^d(\V_n)$ and $w \in S_d$,
making 
$T^d(\V_n)$ into a $(\U_n, \H_d)$-bimodule.
To write this action of $\H_d$ down in a more familiar
way in terms of generators,
let $\alpha = (\alpha_1,\dots,\alpha_d) \in I_n^d$ be a multi-index 
as in $\S$\ref{scombinatorics}.
Define
${M}_\alpha = v_{\alpha_1}\otimes v_{\alpha_2}\otimes \cdots\otimes
v_{\alpha_d}$,
so that
$\{M_{\alpha}\:|\:\alpha \in I_n^d\}$ is the standard basis
for $T^d(\V_n)$.
Then,\begin{align}\label{mbaact}
M_{\alpha}
H_i = 
\left\{
\begin{array}{ll}
M_{\alpha \cdot s_i}
& \hbox{if $\alpha_i < \alpha_{i+1}$,}\\
q^{-1} M_{\alpha}
&\hbox{if $\alpha_i = \alpha_{i+1}$,}\\
M_{\alpha \cdot s_i} - (q-q^{-1}) M_{\alpha}
&\hbox{if $\alpha_i > \alpha_{i+1},$}
\end{array}
\right.
%\\
%M_{\alpha}^*
%H_i = 
%\left\{
%\begin{array}{ll}
%M^*_{\alpha \cdot s_i}
%& \hbox{if $\alpha_i < \alpha_{i+1}$,}\\
%q M_{\alpha}
%&\hbox{if $\alpha_i = \alpha_{i+1}$,}\\
%M^*_{\alpha \cdot s_i} + (q-q^{-1}) M^*_{\alpha}
%&\hbox{if $\alpha_i > \alpha_{i+1}$,}
%\end{array}
%\right.
\end{align}
for each $\alpha \in I_n^d$ and $i=1,\dots,d-1$.
We will also often work with the elements
\begin{equation}\label{bwds}
M_\alpha^* = v_{\alpha_d}\otimes\cdots\otimes v_{\alpha_2}\otimes v_{\alpha_1}
= M_{\alpha\cdot w_d},
\end{equation}
so $\{M_\alpha^*\:|\:\alpha \in I_n^d\}$ is the same basis as before but parametrized in the opposite way.

Now fix a weight $\lambda \in \Lambda_n$ with $|\lambda| = d$ and consider
the $\lambda$-weight space $T_\lambda^d(\V_n)$ of $T^d(\V_n)$.
It is well known, and easy to prove using (\ref{mbaact}) and 
the definition (\ref{permdef}), 
that
the map
\begin{equation*}
\psi_\lambda:T_\lambda^d(\V_n) \rightarrow \M^\lambda, \qquad
M_\alpha \mapsto M_{d(\alpha)},\quad
M^*_\alpha \mapsto M_{w_\lambda d(\alpha)w_d}
\end{equation*}
is an isomorphism of $\H_d$-modules.
The key observation is that the restriction of the 
bar involution on $T^d(\V_n)$ to its $\lambda$-weight space
agrees with the bar involution on $\M^\lambda$ 
under the isomorphism $\psi_\lambda$, i.e. for any $v \in T_\lambda^d(\V_n)$ we have that
$\psi_\lambda(\overline{v}) = \overline{\psi_\lambda(v)}$.
To see this, just note that if $\alpha \in I_\lambda$
has $\alpha_1 \leq \cdots \leq \alpha_d$, 
then $M_\alpha$ is bar invariant by
weight considerations. Since 
$\psi_\lambda(M_\alpha) = M_1$  generates $\M^\lambda$
as an $\H_d$-module and $M_1$ is bar invariant too, it just remains to observe
by (\ref{compr})  that $\overline{v h} =
\overline{v}\,\overline{h}$ for any $v \in T_{\lambda}^d(\V_n)$ and $h \in \H_d$.
We deduce comparing with the opening paragraph of the section
that for every $\alpha \in I_\lambda$
there exist unique bar invariant elements
$L_\alpha$ and $L_\alpha^*$
 in $T_\lambda^d(\V_n)$ such that
$$
{L}_\alpha \in {M}_\alpha + \sum_{\beta \in I_\lambda} q^{-1}\Z[q^{-1}] {M}_\beta,
\quad
L^*_\alpha \in M^*_\alpha + \sum_{\beta \in I_\lambda} q\Z[q] M^*_\beta.
$$
Moreover, 
$\psi_\lambda(L_\alpha) = \underbar{\tilde M}_{d(\alpha)}$
and 
$\psi_\lambda(L_\alpha^*) = \underbar{M}_{w_\lambda d(\alpha)w_d}$.
Let $l_{\alpha,\beta}(q) \in \Z[q^{-1}]$ 
and $l^*_{\alpha,\beta}(q) \in \Z[q]$
denote the coefficients defined from
\begin{align}\label{ekl2}
L_\beta = \sum_{\alpha \in I_\lambda} l_{\alpha,\beta}(q) {M}_\alpha,
\qquad
&L^*_\beta = \sum_{\alpha \in I_\lambda}  l^*_{\alpha,\beta}(q) M^*_\alpha.
\end{align}
These are the same as the coefficients in (\ref{ekl1}),
taking $x = d(\alpha), y = d(\beta)$ and $x = w_\lambda d(\alpha) w_d,
y = w_\lambda d(\beta) w_d$, respectively.
We have now constructed two new bases
$\{L_\alpha\:|\:\alpha \in I_n^d\}$ 
and $\{L^*_\alpha\:|\:\alpha \in I_n^d\}$
for the tensor space $T^d(\V_n)$, which we call the {\em dual canonical} 
and the {\em canonical}
{bases}. In Kashiwara's language, they are
upper and lower global crystal bases, respectively.

Let us recall from \cite{KaG} 
the precise meaning of the previous sentence.
Denote the $\Z[q,q^{-1}]$-submodule of $\V_n$ spanned by 
$v_1,\dots,v_n$ by $\VV_n$; it 
is invariant under the action of
Lusztig's integral form $\UU_n$ for $\U_n$,
i.e. the $\Z[q,q^{-1}]$-subalgebra of $\U_n$
generated by
all $E_i^{(r)}= E_i^r / [r]!, F_i^{(r)}= F_i^r / [r]!$,  $K_i^{\pm 1}$ 
and $\left[ \substack{K_i\\r}\right]= \prod_{s=1}^r \frac{K_i q^{1-s} - K_i^{-1} q^{s-1}}{q^s-q^{-s}}$.
Taking tensor products over $\Z[q,q^{-1}]$, we obtain the
$\Z[q,q^{-1}]$-lattice $T^d(\VV_n)$ in $T^d(\V_n)$.
Next, let $\A_0$ resp. $\A_\infty$ 
be the subring of $\Q(q)$ consisting of all rational functions having
no pole at $q=0$ resp. $q = \infty$, so
$\A_\infty = \overline{\A_0}$.
Let $T^d(\V_n)_0$ resp. $T^d(\V_n)_\infty$ be the 
$\A_0$- resp. $\A_\infty$-submodule of $T^d(\V_n)$ generated
by the elements $\{M_\alpha^*\:|\:\alpha \in I_n^d\}$ resp.
$\{M_\alpha\:|\:\alpha \in I_n^d\}$.
Then, by Kashiwara's tensor product rules \cite{Ka0,Ka1},
$T^d(\V_n)_0$ resp. $T^d(\V_n)_\infty$ is a 
lower resp. upper crystal lattice at $q=0$ resp. $q=\infty$,
and the image of the basis $\{M_\alpha^*\:|\:\alpha \in I_n^d\}$
resp. $\{M_\alpha\:|\:\alpha \in I_n^d\}$ 
in $T^d(\V_n)_0 / q T^d(\V_n)_0$ resp. $T^d(\V_n)_\infty / q^{-1} T^d(\V_n)_\infty$
is a lower resp. upper crystal base at $q=0$ resp. $q=\infty$.
The actions of the lower resp. upper crystal operators
on these crystal bases is described 
by the crystal
$(I_n^d, \tilde e_i, \tilde f_i, \eps_i, \phi_i, \theta)$ 
from $\S$\ref{scombinatorics}.
Finally, the lower resp. upper global crystal base
$\{L_\alpha^*\:|\:\alpha \in I_n^d\}$ resp.
$\{L_\alpha\:|\:\alpha \in I_n^d\}$ is the unique lift of this
local crystal base arising 
from the balanced triple
$(\Q \otimes_{\Z} T^d(\VV_n), T^d(\V_n)_0, \overline{T^d(\V_n)_0})$
resp.
$(\Q \otimes_{\Z} T^d(\VV_n), \overline{T^d(\V_n)_\infty}, {T^d(\V_n)_\infty})$.

The only other thing we want to do in this section is to reprove
the inversion formula for parabolic Kazhdan-Lusztig polynomials
due originally to Douglass \cite{Dou} (see also 
\cite[3.9]{soergel}) in terms of tensor space.
The argument involves an important bilinear form
$(.,.)$ on $T^d(\V_n)$ defined 
by setting $( M_\alpha, \overline{M^*_\beta} ) = \delta_{\alpha,\beta}$
for each $\alpha,\beta \in I_n^d$.
Recall that $\tau:\U_n \rightarrow \U_n$ is the antiautomorphism
with $\tau(E_i) = F_i, \tau(F_i) = E_i$ and $\tau(K_i) = K_i$.
Also let $\tau:\H_d \rightarrow \H_d$ be the antiautomorphism with
$\tau(H_i) = H_{d-i}$.

\begin{lem}\label{mis}
The bilinear form $(.,.)$ is symmetric and 
$( u v h, w ) = ( v, \tau(u) w \tau(h))$ for all
$u \in \U_n$, $h \in \H_d$ and $v,w \in T^d(\V_n)$.
\end{lem}

\begin{proof}
The second part is a routine direct check on generators.
For the first part, we need to show that 
$( \overline{M_\beta^*}, M_\alpha) = \delta_{\alpha,\beta}$.
By (\ref{eq18}), $\overline{M_\beta^*} = 
q^{-\sum_{i < j} (\eps_{\beta_i},\eps_{\beta_j})}
M_\beta H_{w_d}^{-1}$.
Since $\tau(H_{w_d}) = H_{w_d}$, we get that
\begin{align*}
( \overline{M_\beta^*}, M_\alpha)
&= 
q^{-\sum_{i<j} (\eps_{\beta_i},\eps_{\beta_j})}
(  M_\beta H_{w_d}^{-1}, M_\alpha )=q^{-\sum_{i<j} (\eps_{\beta_i},\eps_{\beta_j})}
(  M_\beta, M_\alpha H_{w_d}^{-1} )\\
&= 
q^{\sum_{i<j} ((\eps_{\alpha_i},\eps_{\alpha_j})-(\eps_{\beta_i},\eps_{\beta_j}))}
(  M_\beta, \overline{M_\alpha^*} ) = \delta_{\alpha,\beta}.
\end{align*}
\end{proof}

\begin{thm}\label{due}
$( L_\alpha, L_\beta^* ) = \delta_{\alpha,\beta}$.
\end{thm}

\begin{proof}
Since $L_\beta^*$ is bar invariant, we have by
(\ref{ekl2}) that 
$L_\alpha = \sum_{\gamma} 
l_{\gamma,\alpha}(q) M_\gamma$, 
$L^*_\beta = \sum_{\delta} l^*_{\delta,\beta}(q^{-1}) 
\overline{M^*_\delta}$.
Similarly, 
$L_\beta^* = \sum_{\gamma} l_{\gamma,\beta}^*(q) M_{\gamma\cdot w_d}$,
$L_\alpha = \sum_{\delta} l_{\delta,\alpha}(q^{-1}) \overline{M_{\delta \cdot w_d}^*}$. Hence, by definition of the form,
\begin{align*}
( L_\alpha,L_\beta^* )
= 
\sum_{\gamma} l_{\gamma,\alpha}(q) l_{\gamma,\beta}^*(q^{-1})
\equiv \delta_{\alpha,\beta} &\pmod{q^{-1} \Z[q^{-1}]},\\
( L_\beta^*,L_\alpha )
=
\sum_{\gamma} l_{\gamma,\beta}^*(q) l_{\gamma,\alpha}(q^{-1})
\equiv
\delta_{\alpha,\beta}&\pmod{q\Z[q]}.
\end{align*}
By Lemma~\ref{mis}, we know that
$( L_\alpha, L_\beta^* ) = ( L_\beta^*,L_\alpha )$, so
these two congruences together imply that
$( L_\alpha,L_\beta^*) = \delta_{\alpha,\beta}$.
\end{proof}

\begin{cor}
$\displaystyle
\sum_{\gamma \in I_n^d} l_{\alpha,\gamma}(q) l_{\beta,\gamma}^*(q^{-1}) = \delta_{\alpha,\beta}$.
\end{cor}

\begin{remark}\rm
Let us explain the essential difference between the exposition here and that of 
\cite{FKK}. In that paper, there are two different 
$\U_n$-module structures and two different bar involutions on the 
underlying vector space $T^d(\V_n)$. One of these is 
used to define the dual canonical basis, exactly as here.
The other $\U_n$-module structure, which may be denoted $\widetilde T^d(\V_n)$,
is defined using the comultiplication $\widetilde\Delta$ from $\S$\ref{srmatrices}, and
its compatible bar involution is defined
using 
the corresponding quasi-$R$-matrix $\widetilde\Theta$.
Letting $\widetilde M_\alpha = v_{\alpha_1}\otimes\cdots\otimes v_{\alpha_d}
\in \widetilde T^d(\V_n)$, the canonical basis 
element $\widetilde L_\alpha$ is then 
the unique bar invariant element
lying in $\widetilde M_\alpha + \sum_{\beta \in I_n^d} q^{-1} \Z[q^{-1}] 
\widetilde M_\beta$.
To translate between this and the approach 
followed here, we note that there is a $\U_n$-module isomorphism
$\widetilde T^d(\V_n) \rightarrow T^d(\V_n), \widetilde M_\alpha \mapsto
\overline {M_\alpha^*},
\widetilde L_\alpha\mapsto L_\alpha^*$.
\end{remark}

\section{Symmetric and exterior powers}\label{spowers}

In this section, we define canonical and dual canonical bases
in tensor products of symmetric and exterior powers of $\V_n$,
generalizing the canonical and dual canonical bases of tensor space
from the previous section.
We start with symmetric powers, then summarize the 
necessary changes for exterior powers at the end of the section.
By definition,
the {\em quantum symmetric algebra} 
$S(\V_n)$ is the quotient
of $T(\V_n)$ by the two-sided ideal $I$ generated by the elements
\begin{equation}\label{symdef}
\{v_j \otimes v_i - q^{-1} v_i \otimes v_j\:|\: 1 \leq i < j \leq n\}.
\end{equation}
Clearly,
$I = 
\bigoplus_{d \geq 0} I_d$ where
$I_d = I \cap T^d(\V_n)$.
The {\em $d$th symmetric power} 
$S^d(\V_n)$ is the $d$th homogeneous 
component $T^d(\V_n) / I_d$
of $S(\V_n)$, so 
$S(\V_n) = \bigoplus_{d \geq 0} S^d(\V_n)$.
One checks that 
$I_2$ 
is invariant both under the action of $\U_n$ 
and under the bar involution on $T^2(\V_n)$.
Since $I_2$ generates $I$, it follows that {\em all} $I_d$ 
are invariant under $\U_n$ and under the bar involution. 
Hence, $S^d(\V_n)$ is a $\U_n$-module quotient
of $T^d(\V_n)$, and the bar involution on $T^d(\V_n)$ 
descends to give a compatible bar involution on $S^d(\V_n)$.

We also need the dual object, the {\em $d$th divided power}
$\widetilde{S}^d(\V_n)$. To define this, let
\begin{equation*}
X_d = \sum_{w \in S_d}  q^{\ell(w_d)-\ell(w)} H_w \in \H_d.
\end{equation*}
Then, by definition, $\widetilde{S}^d(\V_n)$ is the $\U_n$-submodule
$T^d(\V_n) X_d$ of $T^d(\V_n)$.
It is well known that $X_d$ is bar invariant, hence
the bar involution on $T^d(\V_n)$ restricts to a well-defined
compatible bar involution on $\widetilde{S}^d(\V_n)$, and also
\begin{equation}\label{hxh}
H_i X_d = q^{-1} X_d = X_d H_i
\end{equation}
for all $i$.
Let $\iota:\widetilde{S}^d(\V_n) \hookrightarrow T^d(\V_n)$ be the inclusion
and $\pi:T^d(\V_n) \twoheadrightarrow S^d(\V_n)$ be the quotient map.
We claim that the bilinear form $(.,.)$ on $T^d(\V_n)$ induces
a well defined pairing $(.,.): {S}^d(\V_n) \times 
\widetilde{S}^d(\V_n) \rightarrow \Q(q)$ with $( \pi(v),w ) = ( v, \iota(w) )$
for all $v \in T^d(\V_n)$ and $w \in \widetilde{S}^d(\V_n)$.
To prove this, we need to show that $( \ker\pi, \im \iota ) = 0$.
By (\ref{symdef}) and (\ref{mbaact}),
$\ker \pi$ is spanned by vectors of the form
$v (H_i-q^{-1})$, while
$\im \iota$ is spanned by vectors of the form
$w X_d$.
Now Lemma~\ref{mis} and (\ref{hxh}) show that
$( v (H_i-q^{-1}),wX_d ) =
( v,wX_d (H_{d-i}-q^{-1}) ) = 0$
proving the claim.

To define the standard bases for the spaces
${S}^d(\V_n)$ and $\widetilde{S}^d(\V_n)$, take
$\alpha \in I_n^d$ with $\alpha_1 \leq \cdots \leq \alpha_d$.
Define $X_\alpha$ to be the bar invariant element 
$\pi(M_\alpha)$ of $S^d(\V_n)$. 
Also, letting $\lambda$ denote the weight $\theta(\alpha)$, set
\begin{equation}
X^*_\alpha = \sum_{\beta \sim\alpha} q^{\ell(\alpha,\beta)}M_{\beta}^*
= \frac{1}{[\lambda_1]!\cdots [\lambda_n]!} M_\alpha X_d,
\end{equation}
where 
$\ell(\alpha,\beta)$ 
denotes the length of the shortest element $w \in S_d$ with
$\beta = \alpha \cdot w$.
Note $X^*_\alpha$ belongs to $\widetilde{S}^d(\V_n)$ and it is bar invariant
(indeed, it coincides with the canonical basis element $L^*_\alpha$).
Using (\ref{symdef}) and (\ref{hxh}) one checks easily that
the vectors $X_\alpha$ and $X_\beta^*$ for all 
weakly increasing $\alpha,\beta \in I_n^d$ span
$S^d(\V_n)$ and $\widetilde{S}^d(\V_n)$, respectively.
Finally, we have that
$$
( X_\alpha, \overline{X_\beta^*} )
= 
\sum_{\gamma\sim\beta} q^{-\ell(\beta,\gamma)}
( M_\alpha, \overline{M_{\gamma}^*} )
= \delta_{\alpha,\beta}.
$$ 
This gives the linear independence needed to show that
$\{X_\alpha\:|\:\alpha \in I_n^d, \alpha_1\leq\cdots\leq \alpha_d\}$
is a basis for $S^d(\V_n)$ and
$\{X^*_\alpha\:|\:\alpha \in I_n^d, \alpha_1\leq\cdots\leq \alpha_d\}$
is a basis for $\widetilde{S}^d(\V_n)$.

Suppose more generally that $\mu \in \Lambda_m$ and $|\mu| = d$.
Consider the $\U_n$-modules
\begin{align*}
{S}^\mu(\V_n) = 
{S}^{\mu_m}(\V_n) \otimes\cdots\otimes {S}^{\mu_1}(\V_n),\qquad\quad
\widetilde{S}^\mu(\V_n) = \widetilde{S}^{\mu_1}(\V_n)\otimes\cdots\otimes
\widetilde{S}^{\mu_m}(\V_n).
\end{align*}
Note
${S}^\mu(\V_n)$ is
a quotient of $T^d(\V_n)$; we let $\pi:T^d(\V_n)\twoheadrightarrow{S}^\mu(\V_n)$
be the quotient homomorphism.
Also, $\widetilde{S}^\mu(\V_n)$ is a submodule of
$T^d(\V_n)$; we let
$\iota:\widetilde{S}^\mu(\V_n) \hookrightarrow T^d(\V_n)$
be the inclusion.
Since all ${S}^{\mu_i}(\V_n)$ and
$\widetilde{S}^{\mu_j}(\V_n)$ possess compatible bar involutions,
we get induced compatible bar involutions on
${S}^\mu(\V_n)$ and $\widetilde{S}^\mu(\V_n)$ by the general construction
explained in $\S$\ref{srmatrices}.
It is immediate from this construction that these bar involutions are
consistent with the one on $T^d(\V_n)$ itself, i.e. the maps $\pi$ and $\iota$
commute with the bar involutions.
As before, the symmetric bilinear form $(.,.)$ on $T^d(\V_n)$
induces a well-defined pairing 
$(.,.):{S}^\mu(\V_n) \times
\widetilde{S}^\mu(\V_n) \rightarrow \Q(q)$ with
$( \pi(v), w ) = ( v, \iota(w))$
for all $v \in T^d(\V_n)$ and $w \in \widetilde{S}^\mu(\V_n)$.

For each $\nu \in \Lambda_n$ with $|\nu| = d$,
there are natural monomial bases
for the $\nu$-weight spaces ${S}_\nu^\mu(\V_n)$
and $\widetilde{S}_\nu^\mu(\V_n)$ of
${S}^\mu(\V_n)$ and $\widetilde{S}^\mu(\V_n)$,
parametrized by the set $\Row(\mu,\nu)$
of row standard tableaux of row shape $\mu$ and weight $\nu$
from $\S$\ref{scombinatorics}.
To write these down, recall 
that $\rho(A)$ is the 
row reading of the tableau $A$. 
For $A \in \Row(\mu,\nu)$,
let $M_A = \pi(M_{\rho(A)})
\in {S}^\mu(\V_n)$ and let
$M_A^* \in \widetilde{S}^\mu(\V_n)$ be the
unique element with
\begin{equation}
\iota(M^*_A) = \sum_{B \sim_\row A} q^{\ell(A,B)} M^*_{\rho(B)},
\end{equation}
writing $\ell(A,B)$ for the 
minimal number of transpositions of neighbouring entries in the same
row needed to get $B$ from $A$.
Then, the vectors
$\{M_A\:|\:A \in \Row(\mu,\nu)\}$ give a basis for
${S}_\nu^\mu(\V_n)$
and the vectors $\{M_A^*\:|\:A \in \Row(\mu,\nu)\}$
give a basis for $\widetilde{S}^\mu_\nu(\V_n)$.
Moreover, the pairing $(.,.)$ satisfies
$( M_A, \overline{M_B^*} ) = \delta_{A,B}$.

We can now introduce the canonical and dual canonical bases.
Recalling the second equivalent definition of the Bruhat ordering on 
$\Row(\mu,\nu)$ from $\S$\ref{scombinatorics}, 
one checks by weight considerations that 
the bar involutions on $S^\mu(\V_n)$, $\widetilde{S}^\mu(\V_n)$
satisfy
\begin{align*}
\overline{M_A} &= M_A + \text{(a $\Z[q,q^{-1}]$-linear combination
of $M_B$'s for $B < A$)},\\
\overline{M_A^*} &= M_A^* + \text{(a $\Z[q,q^{-1}]$-linear combination
of $M_B^*$'s for $B > A$)}.
\end{align*}
Hence by \cite[1.2]{Du}, we deduce that for every $A \in \Row(\mu,\nu)$ there are unique bar invariant elements
$L_A \in {S}^\mu(\V_n)$ and $L_A^* \in \widetilde{S}^\mu(\V_n)$ 
such that
$$
L_A \in M_A + \sum_{B \in \Row(\mu,\nu)} q^{-1}\Z[q^{-1}] M_B,
\quad
{L}^*_A \in {M}^*_A + \sum_{B \in \Row(\mu,\nu)} q\Z[q] {M}^*_B.
$$
As before, we introduce notation for the coefficients:
\begin{align}\label{ekl3}
L_B = \sum_{A \in \Row(\mu,\nu)}  l_{A,B}(q) M_A,\qquad
&L^*_B = \sum_{A \in \Row(\mu,\nu)} l^*_{A,B}(q) {M}^*_A.
\end{align}
The polynomials $l_{A,B}(q) \in \Z[q^{-1}]$, $l^*_{A,B}(q) \in \Z[q]$
satisfy $l_{A,A}(q) = l_{A,A}^*(q) = 1$ and $l_{A,B}(q) = 0$
unless $A \leq B$, $l_{A,B}^*(q) = 0$ unless $A \geq B$.
We have now constructed two new bases
$\{L_A\:|\:A \in \Row(\mu,\nu)\}$ 
for ${S}_\nu^\mu(\V_n)$ and $\{L_A^*\:|\:A \in \Row(\mu,\nu)\}$
for $\widetilde{S}_\nu^\mu(\V_n)$, which we 
call the {\em dual canonical} and the {\em canonical}
{bases}, respectively.
They are upper and lower global crystal bases in the sense of \cite{KaG},
the precise meaning of this phrase being just like 
in the previous section.
We just note that 
the constructions just described can be carried out equally well
over the ring $\Z[q,q^{-1}]$, to obtain the natural integral forms
$S^\mu(\VV_n)$ and $\widetilde S^\mu(\VV_n)$.
Thus, $S^\mu(\VV_n)$ is the free $\Z[q,q^{-1}]$-module 
with basis given either by the $M_A$'s or by the $L_A$'s, 
$\widetilde S^\mu(\VV_n)$ is the
free $\Z[q,q^{-1}]$-module with basis given either by the
$M_A^*$'s or by the $L_A^*$'s.
Both are invariant under the action of Lusztig's $\Z[q,q^{-1}]$-form
$\UU_n$.

\begin{thm}\label{bermuda}
For $A \in \Row(\mu,\nu)$ we have that
$\iota(L^*_A) = L^*_{\rho(A)}$
and $\pi(L_{\rho(A)}) = L_A$.
Moreover, if $\alpha \in I_\nu$ is not equal to
$\rho(A)$ for any $A \in \Row(\mu,\nu)$, then
$\pi(L_\alpha) = 0$.
Hence, $( L_A, L_B^* ) = \delta_{A,B}$
for all $A, B \in \Row(\mu,\nu)$.
\end{thm}

\begin{proof}
Note for $A \in \Row(\mu,\nu)$ that $\iota(L_A^*)$ 
is bar invariant and it equals $M^*_{\rho(A)}$ plus
a $q\Z[q]$-linear combination of 
$M^*_{\beta}$'s. Hence, $\iota(L^*_A) = L^*_{\rho(A)}$.
Similarly, $\pi(L_{\rho(A)})$ is bar invariant and it equals
$M_A$ plus a $q^{-1}\Z[q^{-1}]$-linear combination of 
$M_B$'s. Hence, it equals $L_A$.
Moreover, if $\alpha \in I_\nu$ is not equal to $\rho(A)$ for any
$A \in \Row(\mu,\nu)$, then 
$\pi(L_{\alpha})$ is bar invariant and it is a 
$q^{-1}\Z[q^{-1}]$-linear combination of $M_B$'s. Hence, it must be zero.
Finally, for any $A, B \in \Row(\mu,\nu)$,
we get that $( L_A, L_B^* ) = ( L_{\rho(A)}, 
L_{\rho(B)}^* ) = \delta_{A,B},$ 
using Theorem~\ref{due}.
\end{proof}

\begin{cor}\label{if}
$\displaystyle\sum_{C \in \Row(\mu,\nu)} l_{A, C}(q) l_{B,C}^*(q^{-1}) = \delta_{A, B}$.
\end{cor}

\begin{cor}\label{nc1}
$\,\displaystyle 
l_{A,B}(q) = \sum_{C \sim_\row A} q^{-\ell(A,C)} l_{\rho(C), \rho(B)}(q),\quad
l_{A,B}^*(q) = l^*_{\rho(A), \rho(B)}(q)$.
\end{cor}

\begin{remark}\label{okl}\rm
Using Corollary~\ref{nc1}, the identification of the polynomials in (\ref{ekl1})
and (\ref{ekl2}), and \cite[2.6, 3.4]{soergel},
we obtain the following formulae relating the polynomials $l_{A,B}(q)$
and $l_{A,B}^*(q)$ directly to the original Kazhdan-Lusztig polynomials
$P_{x,y}(t) \in \Z[t]$ from \cite{KL}: 
\begin{align*}
l_{A,B}(q) &= q^{\ell(y)-\ell(x)} \sum_{z \in S_\nu x S_\mu}
(-1)^{\ell(z)+\ell(y)}
P_{zw_d,yw_d}(q^2),\\
l^*_{A,B}(q) &= q^{\ell(y)-\ell(x)} P_{x,y}(q^{-2}),\label{lf2}
\end{align*}
where $x = d(\rho(A))w_d$ and $y=d(\rho(B))w_d$.
Using these formulae, one sees that
the inversion formula from Corollary~\ref{if} 
is the same as \cite[1.3]{Du2}.
\end{remark}

We now turn our attention to exterior powers. The proofs 
are all the same as the above proofs for symmetric powers, so we omit them.
Note however that it is necessary throughout to 
interchange the roles of canonical and dual canonical bases.
The {\em quantum exterior algebra}
$\widetilde{\bw}(\V_n)$ is the 
quotient of $T(\V_n)$ by the homogeneous two-sided
ideal $J = \bigoplus_{d \geq 0} J_d$ generated by the elements
\begin{equation*}
\{v_j \otimes v_i + q v_i \otimes v_j\:|\:1 \leq i < j \leq n\} \cup \{v_i \otimes v_i\:|\:1 \leq i \leq  n\}.
\end{equation*}
Let $\widetilde{\bw}^d(\V_n)$ be the
$d$th homogeneous component $T^d(\V_n) / J_d$.
It is a $\U_n$-module and inherits a compatible bar involution from the
one on $T^d(\V_n)$. Although one usually calls $\widetilde{\bw}^d(\V_n)$ the
$d$th exterior power, we prefer here to reserve 
that name for the (isomorphic) dual object
${\bw}^d(\V_n) = T^d(\V_n) Y_d$
where
\begin{equation*}
Y_d = \sum_{w \in S_d} (-q)^{\ell(w)-\ell(w_d)} H_w \in \H_d.
\end{equation*}
Since $Y_d$ is bar invariant, the bar involution 
on $T^d(\V_n)$ restricts to a compatible bar involution on ${\bw}^d(\V_n)$.
Recall also that
$H_i Y_d = -q Y_d = Y_d H_i$
for all $i=1,\dots,d-1$.
For $\alpha \in I_n^d$ with $\alpha_1 > \cdots > \alpha_d$, 
let
\begin{equation}\label{ya}
Y_\alpha = \sum_{\beta\sim\alpha} (-q)^{-\ell(\alpha,\beta)}
M_\beta = M_\alpha^* Y_d \in \bw^d(\V_n).
\end{equation}
Also let $Y_\alpha^*$ be the image of $M_\alpha^*$ in the quotient
$\widetilde\bw^d(\V_n)$,
Both $Y_\alpha$ and $Y_\alpha^*$ are bar invariant (indeed $Y_\alpha = 
L_\alpha$).
The vectors 
$\{Y_\alpha\:|\:\alpha \in I_n^d, \alpha_1 > \cdots > \alpha_d\}$ give
a basis for ${\bw}^d(\V_n)$ and the vectors
$\{Y^*_\alpha\:|\:\alpha \in I_n^d, \alpha_1 > \cdots > \alpha_d\}$ give
a basis for $\widetilde{\bw}^d(\V_n)$. 

Now take $\mu \in \Lambda_l$ and $\nu \in \Lambda_n$
with $|\mu|=|\nu| = d$.
Consider the $\U_n$-modules
$$
{\bw}^\mu(\V_n) = {\bw}^{\mu_1}(\V_n)\otimes\cdots\otimes
{\bw}^{\mu_l}(\V_n),\quad\quad
\widetilde{\bw}^\mu(\V_n) = \widetilde{\bw}^{\mu_l}(\V_n)\otimes\cdots\otimes
\widetilde{\bw}^{\mu_1}(\V_n).
$$
We write ${\bw}^\mu_\nu(\V_n)$ and $\widetilde{\bw}^\mu_\nu(\V_n)$ for the
$\nu$-weight spaces of these modules.
Also let $\iota:{\bw}^\mu(\V_n) \hookrightarrow T^d(\V_n)$ be the natural
inclusion and $\pi:T^d(\V_n) \twoheadrightarrow \widetilde{\bw}^\mu(\V_n)$ be the natural quotient map. 
There are compatible bar involutions on ${\bw}^\mu(\V_n)$ and on
$\widetilde{\bw}^\mu(\V_n)$, consistent with the bar involution
on $T^d(\V_n)$, and the form $(.,.)$ induces a pairing
$(.,.):\bw^\mu(\V_n) \times \widetilde{\bw}^\mu(\V_n)
\rightarrow \Q(q)$.
To define bases here, recall the definition of the set
$\Col(\mu,\nu)$ and the column reading $\gamma(A)$
of a tableau of column shape $\mu$ from $\S$\ref{scombinatorics}.
For $A \in \Col(\mu,\nu)$, 
let $N_A$ denote the unique element of $\bw^\mu(\V_n)$ with
\begin{equation}\label{nadef}
\iota(N_A) = \sum_{B \sim_\col A} (-q)^{-\ell'(A,B)} M_{\gamma(B)},
\end{equation}
where $\ell'(A,B)$ denotes $\ell(A',B')$.
Let $N_A^* = \pi(M_{\gamma(A)}^*)$.
Then, $\{N_A\:|\:A \in \Col(\mu,\nu)\}$ is a basis for 
${\bw}^\mu_\nu(\V_n)$ and $\{N_A^*\:|\:A \in \Col(\mu,\nu)\}$ is a basis
for $\widetilde{\bw}_\nu^\mu(\V_n)$. 
Moreover, the pairing $(.,.)$
satisfies
$( N_A, \overline{N_B^*} ) = \delta_{A,B}$.
We have that
\begin{align*}
\overline{N_A} &= N_A + \text{(a $\Z[q,q^{-1}]$-linear combination
of $N_B$'s for $B <' A$)},\\
\overline{N_A^*} &= N_A^* + \text{(a $\Z[q,q^{-1}]$-linear combination
of $N_B^*$'s for $B >' A$)}.
\end{align*}
Hence by \cite[1.2]{Du} there are unique bar invariant elements
$K_A \in \bw^\mu(\V_n)$ and $K_A^* \in \widetilde{\bw}^\mu(\V_n)$ 
for each $A \in \Col(\mu,\nu)$
such that
$$
K_A \in N_A + \sum_{B \in \Col(\mu,\nu)} q^{-1}\Z[q^{-1}] N_B,
\qquad
{K}^*_A \in {N}^*_A + \sum_{B \in \Col(\mu,\nu)} q\Z[q] {N}^*_B.
$$
We let
\begin{align}\label{ekl4}
K_B = \sum_{A \in \Col(\mu,\nu)}  k_{A,B}(q) N_A,\qquad
&K^*_B = \sum_{A \in \Col(\mu,\nu)} k^*_{A,B}(q) {N}^*_A.
\end{align}
Note $k_{A,B}(q) \in \Z[q^{-1}]$ and $k_{A,B}^*(q) \in \Z[q]$
satisfy $k_{A,A}(q) = k_{A,A}^*(q) = 1$ and $k_{A,B}(q) = 0$
unless $A \leq' B$, $k_{A,B}^*(q) = 0$ unless $A \geq' B$.
We have now constructed 
bases
$\{K_A\:|\:A \in \Col(\mu,\nu)\}$ 
for $\bw_\nu^\mu(\V_n)$
and $\{K_A^*\:|\:A \in \Col(\mu,\nu)\}$
for $\widetilde{\bw}_\nu^\mu(\V_n)$,
which are the {\em dual canonical} ($=$ upper global crystal)
and {\em canonical} ($=$ lower global crystal) bases, respectively.
Finally, we note that the $\Z[q,q^{-1}]$-submodules of
$\bw^\mu(\V_n)$ and $\widetilde\bw^\mu(\V_n)$ spanned by these bases
give the natural integral forms
$\bw^\mu(\VV_n)$ and $\widetilde\bw^\mu(\VV_n)$, which are invariant under the action of
$\UU_n$.

\begin{thm}\label{bermuda2}
For $A \in \Col(\mu,\nu)$ we have that
$\iota(K_A) = L_{\gamma(A)}$ and $\pi(L_{\gamma(A)}^*) = K_A^*$.
Moreover, if $\alpha \in I_\nu$ is not equal to $\gamma(A)$
for any $A \in \Col(\mu,\nu)$, then $\pi(L_\alpha^*)=0$.
Hence, $( K_A, K_B^* ) = \delta_{A,B}$
for all $A, B \in \Col(\mu,\nu)$.
\end{thm}

\begin{cor}\label{kinv}
$\displaystyle\sum_{C \in \Col(\mu,\nu)} k_{A, C}(q) k_{B,C}^*(q^{-1}) = \delta_{A, B}$.
\end{cor}

\begin{cor}\label{nc2}
$\displaystyle k_{A,B}(q) = l_{\gamma(A), \gamma(B)}(q)$,
$\displaystyle k^*_{A,B}(q) = \!\!\sum_{C \sim_\col A} (-q)^{\ell'(A,C)} l^*_{\gamma(C), \gamma(B)}(q)$.
\end{cor}

\begin{remark}\label{dol}\rm
Using
 Corollary~\ref{nc2} and \cite[2.6, 3.4]{soergel}
as before,
we get that
\begin{align*}
k_{A,B}(q) &=
(-q)^{\ell(x)-\ell(y)} \sum_{z \in S_\nu} (-1)^{\ell(z)}
P_{zx,y}(q^2),\\
k_{A,B}^*(q) &= 
(-q)^{\ell(x)-\ell(y)}
\sum_{z \in \D_\nu \cap S_\nu x S_\mu}
(-1)^{\ell(z)+\ell(y)} P_{z w_d, y w_d}(q^{-2}),
\end{align*}
where $x=d(\gamma(A))$ and $y = d(\gamma(B))$.
\end{remark}

\section{The quantized coordinate algebra}\label{sqca}

Consider now the tensor product
$T^m(S(\V_n)) = S(\V_n)\otimes\cdots\otimes S(\V_n)$
of $m$ copies of the symmetric algebra
$S(\V_n)$.
We obviously have that
\begin{equation*}
T^m(S(\V_n)) = \bigoplus_{\substack{(\mu,\nu) \in \Lambda_m\times\Lambda_n \\ |\mu|=|\nu|}} 
S_\nu^{\mu}(\V_n)
\end{equation*}
so $T^m(S(\V_n))$ has standard basis 
$\{M_A\:|\:A \in \bc\Row(\mu,\nu)\}$
and 
dual canonical basis $\{L_A\:|\:A \in \bc\Row(\mu,\nu)\}$,
where throughout the section $\bc$ 
denotes the union over all pairs
$(\mu,\nu) \in \Lambda_m \times \Lambda_n$ with $|\mu| = |\nu|$.
Because $S(\V_n)$ is a 
polynomial $\U_n$-algebra equipped with a compatible bar involution, 
$T^m(S(\V_n))$
also has a canonical algebra structure and a compatible bar involution, 
defined 
as at the end of $\S$\ref{srmatrices}.
It is well known that this algebra coincides with
the quantized coordinate algebra 
$\O_q(\M_{m,n})$ of the variety
$\M_{m,n}$ of $m \times n$ matrices, that is,
the $\Q(q)$-algebra
on generators
$\{x_{i,j}\:|\:i=1,\dots,m, j = 1,\dots,n\}$ subject only to the relations
\begin{align*}
x_{i,j} x_{k,l} &= x_{k,l} x_{i,j} &(i < k, j > l)\\
x_{i,j} x_{k,l} &= x_{k,l} x_{i,j} + (q-q^{-1})  
x_{k,j}x_{i,l}&(i < k, j < l)\\
x_{i,j} x_{k,j} &= q x_{k,j} x_{i,j}&(i < k)\\
x_{i,j} x_{i,l} &= q x_{i,l} x_{i,j}&(j < l)
\end{align*}
for all $1 \leq i,k \leq m$ and $1 \leq j,l \leq n$.
A proof is written down in \cite[4.2]{BZw} (see also \cite{Zhang}), but still
we repeat the argument in Theorem~\ref{qca} 
below since some of our choices are slightly different.
First, we develop a little more combinatorial language.
Recall from $\S$\ref{scombinatorics} that if $|\mu| = |\nu| = d$, then
the set $\Row(\mu,\nu)$ is in canonical bijection with the set
$(I_\mu\times I_\nu) / S_d$.
We call elements $(\alpha,\beta) \in I_\mu \times I_\nu$
{\em double indexes}. For such a double index $(\alpha,\beta)$, 
introduce the monomial
$$
M_{\alpha,\beta} := x_{\alpha_1,\beta_1} x_{\alpha_2,\beta_2}
\cdots x_{\alpha_d,\beta_d} \in \O_q(\M_{m,n}).
$$
We say that a double
index $(\alpha,\beta)$ is {\em initial}
if $\alpha_1 \geq \cdots \geq \alpha_d$ and $\beta_i \leq \beta_{i+1}$
whenever $\alpha_i = \alpha_{i+1}$, and
{\em terminal} 
if $\beta_1 \leq \cdots \leq \beta_d$ and $\alpha_i \geq \alpha_{i+1}$
whenever $\beta_i = \beta_{i+1}$.
Let $(I_\mu \times I_\nu)^+$ and $(I_\mu \times I_\nu)^-$ denote the
sets of all initial and terminal double indexes in $I_\mu\times I_\nu$,
respectively. These give two distinguished choices of representatives
for the orbits in $(I_\mu\times I_\nu)/S_d$.
The canonical bijection $\Row(\mu,\nu) \rightarrow (I_\mu \times I_\nu)^+$  maps
$A \in \Row(\mu,\nu)$ to the unique initial double index
$(\alpha,\beta) \in (I_\mu\times I_\nu)^+$ with $\beta = \rho(A)$;
see the end of the paragraph after (\ref{eg}) for an example.

\begin{thm}\label{qca}
There is an algebra isomorphism
$\psi:T^m(S(\V_n)) \rightarrow \O_q(\M_{m,n})$ 
defined
for 
any
$(\mu,\nu) \in \Lambda_m\times \Lambda_n$ with
$|\mu|=|\nu|$
and $A \in \Row(\mu,\nu)$
 by
$\psi(M_A) = M_{\alpha,\beta}$,
where $(\alpha,\beta) \in (I_\mu\times I_\nu)^+$ is defined from
$\beta = \rho(A)$. 
In particular, the monomials
$\{M_{\alpha,\beta}\:|\:(\alpha,\beta) \in \bc (I_\mu\times I_\nu)^+\}$
form a basis for $\O_q(\M_{m,n})$.
\end{thm}

\begin{proof}
One checks relations using (\ref{eq18}),
(\ref{barinv}) and (\ref{symdef}) to see that there is a well-defined
algebra homomorphism
$\O_q(\M_{m,n}) \rightarrow T^m(S(\V_n)$ mapping
the generator $x_{i,j}$ 
to $1 \otimes \cdots \otimes 1 \otimes v_j \otimes 1 \otimes \cdots \otimes 1$ (where $v_j$ appears in the $(m+1-i)$th tensor position) for each
$1 \leq i \leq m, 1 \leq j \leq n$.
This maps $M_{\alpha,\beta}$ to $M_A$, hence it is an isomorphism
since the vectors $\{M_A\:|\:A \in \bc \Row(\mu,\nu)\}$
are linearly independent and by the relations the monomials
$\{M_{\alpha,\beta}\:|\:(\alpha,\beta) \in \bc (I_\mu\times I_\nu)^+\}$
span $\O_q(\M_{m,n})$.
The map $\psi$ is the inverse isomorphism.
\end{proof}

Let us view
$\O_q(\M_{m,n})$ as an $X_m \times X_n$-graded algebra
by declaring that the generator $x_{i,j}$ is of
degree $(\eps_i,\eps_j) \in X_m \times X_n$.
Thus, $$
\O_q(\M_{m,n}) = \bigoplus_{\substack{(\mu,\nu)\in\Lambda_m\times\Lambda_n\\|\mu|=|\nu|}}
\O_q(\M_{m,n})_{\mu,\nu}
$$
where 
$\O_q(\M_{m,n})_{\mu,\nu}$ has basis
$\{M_{\alpha,\beta}\:|\:(\alpha,\beta) \in (I_\mu\times I_\nu)^+\}$.
From now on, we're going to {\em  identify} 
$\O_q(\M_{m,n})_{\mu,\nu}$ with the $\nu$-weight space
$S^\mu_\nu(\V_n)$ of $S^\mu(\V_n)$ via the isomorphism $\psi$
from Theorem~\ref{qca}. Thus, for $(\alpha,\beta) \in (I_\mu\times I_\nu)^+$,
the monomial $M_{\alpha,\beta}$ is identified with 
$M_A$, where $A \in \Row(\mu,\nu)$ is defined from
$\beta = \rho(A)$.
The next result gives a direct description
of the bar involution on $\O_q(\M_{m,n})$
arising from this identification.

\begin{thm}
\label{qbar} The bar involution on $\O_q(\M_{m,n})$ 
is the unique antilinear map such that
$\overline{x_{i,j}} = 
x_{i,j}$ for all $1 \leq i\leq m, 
1 \leq j\leq n$ and
$\overline{xy} = q^{(\mu,\bar\mu)-(\nu,\bar\nu)} 
\overline{y}\, \overline{x}$
for all
$x \in \O_q(\M_{m,n})_{\mu,\nu}$ and
$y \in \O_q(\M_{m,n})_{\bar\mu,\bar\nu}$.
Moreover, for $(\alpha,\beta) \in (I_\mu \times I_\nu)^+$, we have that
$\overline{M_{\alpha,\beta}} = M_{\alpha',\beta'}$
where $(\alpha', \beta')\in (I_\mu\times I_\nu)^-$ is the unique
terminal double index lying in the same $S_d$-orbit as $(\alpha,\beta)$.
\end{thm}

\begin{proof}
Let $*$ be the twisted multiplication on $T^m(S(\V_n))$
from Lemma~\ref{barprop}. One checks that the twisted multiplication on
$S(\V_n)$ itself satisfies $x * y = q^{d\bar d} xy$ for
$x \in S^d(\V_n)$ and $y \in S^{\bar d}(\V_n)$.
Hence if $x_m\otimes\cdots\otimes x_1 \in S^{\mu_m} 
(\V_n) \otimes\cdots\otimes
S^{\mu_1}(\V_n)$ is of weight $\nu$
and
 $y_m\otimes\cdots\otimes y_1 \in S^{\bar \mu_m} 
(\V_n) \otimes\cdots\otimes
S^{\bar\mu_1}(\V_n)$ is of weight $\bar\nu$, we have that
$$
\overline{(y_m \otimes \cdots \otimes y_1)} *
\overline{(x_m \otimes \cdots \otimes x_1)}
= q^{\mu_1 \bar\mu_1+\cdots+\mu_m\bar\mu_m}
(\overline{y_m \otimes \cdots \otimes y_1})\:
(\overline{x_m \otimes \cdots \otimes x_1}),
$$
so by Lemma~\ref{barprop} 
$$
\overline{(x_m\otimes\cdots\otimes x_1)(y_m\otimes\cdots\otimes y_1)}
= q^{(\mu,\bar\mu)-(\nu,\bar\nu)} 
(\overline{y_m \otimes \cdots \otimes y_1})\:
(\overline{x_m \otimes \cdots \otimes x_1}).
$$
Clearly 
$\overline{x_{i,j}} = x_{i,j}$, so this proves the first statement of the lemma.
The second can then be deduced by induction on $d$
using the defining relations in $\O_q(\M_{m,n})$.
\end{proof}

Using Theorem~\ref{qbar}, we can also give a direct characterization of 
the dual canonical basis
$\{L_A\:|\:A \in \bc \Row(\mu,\nu)\}$ 
of $\O_q(\M_{m,n})$ arising from its identification with $T^m(S(\V_n))$.
We often denote this basis instead by
$\{L_{\alpha,\beta}\:|\:(\alpha,\beta) \in \bc(I_\mu\times I_\nu)^+\}$,
where
for each initial double index $(\alpha,\beta) \in (I_\mu\times I_\nu)^+$,
$L_{\alpha,\beta}$ is the unique bar invariant element of 
$\O_q(\M_{m,n})$ with the property that
\begin{equation}\label{l1}
L_{\alpha,\beta} \in M_{\alpha,\beta} + \sum_{(\alpha',\beta') \in (I_\mu\times I_\nu)^+} q^{-1}\Z[q^{1}] M_{\alpha',\beta'}.
\end{equation}
Applying the bar involution using Theorem~\ref{qbar}, 
it is equally natural to parametrize
this basis by terminal double indexes: it is the
basis $\{L_{\alpha,\beta}\:|\:(\alpha,\beta)\in\bc(I_\mu\times I_\nu)^-\}$
where
$L_{\alpha,\beta}$ is the unique bar invariant element of $\O_q(\M_{m,n})$ with
\begin{equation}\label{l2}
L_{\alpha,\beta} \in M_{\alpha,\beta} + 
\sum_{(\alpha',\beta') \in (I_\mu\times I_\nu)^-} q\Z[q] M_{\alpha',\beta'},
\end{equation}
for $(\alpha,\beta)\in (I_\mu\times I_\nu)^-$.

\begin{remark}
\rm
One finds this elementary approach to the definition of the
dual canonical basis
of $\O_q(\M_{n,n})$ already in work of 
Zhang \cite{Zh}. Actually, Zhang uses an even simpler
modified definition of the
bar involution: his dual canonical basis is invariant instead under
the antilinear algebra antiautomorphism 
$\phi:\O_q(\M_{m,n}) \rightarrow \O_q(\M_{m,n})$
defined by 
$\phi(x_{i,j}) = x_{i,j}$ for all $1 \leq i \leq m, 1 \leq j \leq n$.
This is related to the bar involution defined here by the
equation $\phi(x) = q^{((\nu,\nu)-(\mu,\mu))/2} \overline{x}$
for $x \in \O_q(\M_{m,n})_{\mu,\nu}$. The dual canonical basis in \cite{Zh}
is equal to the dual canonical basis here
up to multiplication by a power of $q$.
\end{remark}

In the remainder of the section, we wish to record proofs of some further 
properties of this dual canonical basis, all of which are 
known but surprisingly hard to find explicitly in the literature.
They were explained to me by Arkady Berenstein, who describes
them as ``folklore''.
First, to compensate for the asymmetry of our identification of
$\O_q(\M_{m,n})$ with $\T^m(S(\V_n))$,
there is an obvious duality
between $m \times n$ matrices and $n \times m$ matrices:
let $\tau:\O_q(\M_{m,n})\rightarrow \O_q(\M_{n,m})$ be the antilinear 
algebra antiisomorphism defined on generators
by $\tau(x_{i,j}) = x_{j,i}$, i.e.
\begin{equation}\label{b1}
\tau(M_{\alpha,\beta}) = M_{\beta\cdot w_d,\alpha \cdot w_d}
\end{equation}
for $\alpha,\beta \in I_n^d$.
Note if $(\alpha,\beta)$ is initial, then
$(\beta\cdot w_d, \alpha\cdot w_d)$ is terminal.
Moreover, 
by Theorem~\ref{qbar}, we have that $\tau(\overline{x}) = \overline{\tau(x)}$
for all $x \in \O_q(\M_{m,n})$. Hence, for $(\alpha,\beta) \in (I_\mu
\times I_\nu)^+$,
$\tau(L_{\alpha,\beta})$ is bar invariant, and the definitions
(\ref{l1})--(\ref{l2}) now imply that
\begin{equation}\label{b2}
\tau(L_{\alpha,\beta}) = L_{\beta\cdot w_d, \alpha\cdot w_d}.
\end{equation}
The equations (\ref{b1})--(\ref{b2}) imply some symmetry
in the transition matrices from (\ref{ekl3}). To write this down,
define a bijection $\tau:\Row(\mu,\nu)\rightarrow 
\Row(\nu,\mu)$ by letting $\tau(A)$
be the unique element of $\Row(\nu,\mu)$
such that the number of entries on the $i$th row of $\tau(A)$
that equal $j$ is the same as the number of entries on the $j$th row
of $A$ that equal $i$, for each $A \in \Row(\mu,\nu)$.

\begin{lem}\label{symmetry}
$l_{A,B}(q) = l_{\tau(A), \tau(B)}(q),\quad
l_{A,B}^*(q) = l_{\tau(A),\tau(B)}^*(q)$.
\end{lem}

\begin{proof}
The first equalty is immediate from 
(\ref{b1})--(\ref{b2}) and the definitions; the second then follows
using the inversion formula from Corollary~\ref{if}.
\end{proof}

Next, we derive a closed
formula for the dual canonical basis of $\O_q(\M_{2,n}) = S(\V_n)\otimes 
S(\V_n)$,
i.e. the dual canonical basis elements
$L_A$ for row standard tableaux
$A$ with just {\em two rows}. Given $r, s \geq 0$
and integers $1 \leq a_1,\dots,a_r,b_1,\dots,b_s \leq n$ 
we will use the shorthand
$M(\substack{a_1\cdots a_r \\ b_1\cdots b_s})$ resp.
$L(\substack{a_1\cdots a_r \\ b_1\cdots b_s})$
for $M_A$ resp. $L_A$,
where $A$ is the row standard tableau with entries 
$a_1,\dots,a_r$ on the top row and $b_1,\dots,b_s$ on 
the bottom row (arranged of course into weakly increasing order).
For example, we have that
$M(\substack{a\\b}) = x_{2,a} x_{1,b}$, and
\begin{equation}\label{lab}
L(\substack{a\\b}) = 
\left\{
\begin{array}{ll}
M(\substack{a\\b}) - q^{-1} M(\substack{b\\a})
&\text{if $a > b$,}\\
M(\substack{a\\b})&\text{if $a \leq b$.}
\end{array}\right.
\end{equation}

\begin{lem}\label{ms}
Let $1 \leq a_1,\dots,a_r,b_1,\dots,b_s,a,b \leq n$ such that
$a > b$, $a_1 \leq \cdots \leq a_r$ and
$b_1 \leq \cdots \leq b_s$.
Assume that
$a_i \notin \{b+1,\dots,a-1\}$ 
for each $i=1,\dots,r$ and
$b_j \notin \{b+1,\dots,a-1\}$
 for each $j=1,\dots,s$.
Then,
\begin{align*}
L(\substack{a_1\cdots a_r a \\ b_1\cdots b_r b})
&=
q^{\#\{i\:|\:a_i > a\}+\#\{j\:|\:b_j >a\}}
L(\substack{a_1\cdots a_r \\ b_1 \cdots b_s}) L(\substack{a\\b})\\&=
q^{\#\{i\:|\:a_i < b\}+\#\{j\:|\:b_j <b\}}
L(\substack{a\\b})L(\substack{a_1\cdots a_r \\ b_1 \cdots b_s}).
\end{align*}
\end{lem}

\begin{proof}
Let $\omega:\O_q(\M_{2,n}) \rightarrow \O_q(\M_{2,n})$
be the linear map defined by
\begin{equation*}
\omega(M(\substack{a_1\cdots a_r \\ b_1\cdots b_s}))
=
M(\substack{a_1\cdots a_r a \\ b_1\cdots b_s b})
- q^{-1-\#\{i\:|\:a_i = a\} - \#\{j\:|\:b_j=b\}}
M(\substack{a_1\cdots a_r b \\ b_1\cdots b_s a})
\end{equation*}
for any $r,s \geq 0$ and
$1 \leq a_1,\dots,a_r,b_1,\dots,b_s \leq n$.
Using the relations and (\ref{lab}),
one checks that
$$
x_{1,b_j} 
L(\substack{a\\b})
=
\left\{
\begin{array}{ll}
q^{-1}
L(\substack{a\\b})
x_{1,b_j} 
&\hbox{if $b_j > a > b$,}\\
L(\substack{a\\b})
x_{1,b_j} 
&\hbox{if $b_j = a > b$,}\\
L(\substack{a\\b})
x_{1,b_j} 
&\hbox{if $a > b=b_j$,}\\
qL(\substack{a\\b})x_{1,b_j} 
&\hbox{if $a > b>b_j$.}
\end{array}
\right.
$$
Hence, recalling that
$M(\substack{a_1\cdots a_r \\ b_1 \cdots b_s}) =
x_{2,a_1}x_{2,a_2}
\cdots x_{2,a_r} x_{1,b_1} x_{1,b_2}\cdots x_{1,b_s}$,
$$
M(\substack{a_1\cdots a_r \\ b_1\cdots b_s})
L(\substack{a\\b}) =
q^{\#\{j\:|\:b_j < b\} - \#\{j\:|\:b_j > a\}}
x_{2,a_1}\cdots x_{2,a_r} 
L(\substack{a\\b})
x_{1,b_1}\cdots x_{1,b_s}
$$
Moreover,
\begin{align*}
x_{2,a_1}&\cdots 
x_{2,a_r} 
L(\substack{a\\b})
x_{1,b_1}\cdots x_{1,b_s} 
=
x_{2,a_1}\cdots x_{2,a_r} 
(x_{2,a}x_{1,b}-q^{-1}x_{2,b}x_{1,a})
x_{1,b_1}\cdots x_{1,b_s}\\
&\:\:\:\:=
q^{-\#\{i\:|\:a_i > a\}-\#\{j\:|\:b_j < b\}}
M(\substack{a_1\cdots a_r a \\ b_1 \cdots b_s b})-q^{-1-\#\{i\:|\:a_i > b\}-\#\{j\:|\:b_j < a\}}
M(\substack{a_1\cdots a_r b \\ b_1 \cdots b_s a}).
\end{align*}
Hence
$\omega(M(\substack{a_1\cdots a_r\\b_1\cdots b_s}))
=
q^{\#\{i\:|\:a_i > a\}+\#\{j\:|\:b_j >a\}}
M(\substack{a_1\cdots a_r\\b_1\cdots b_s})L(\substack{a\\b})$.
A similar argument using instead the relations
$$
L(\substack{a\\b}) x_{2,a_i} =
\left\{
\begin{array}{ll}
qx_{2,a_i} L(\substack{a\\b})
&\hbox{if $a_i > a > b$},\\
x_{2,a_i} L(\substack{a\\b})
&\hbox{if $a_i = a > b$},\\
x_{2,a_i} L(\substack{a\\b})
&\hbox{if $a > b=a_i$},\\
q^{-1}x_{2,a_i} L(\substack{a\\b})
&\hbox{if $a > b>a_i$}
\end{array}\right.
$$
shows that
$\omega(M(\substack{a_1\cdots a_r\\b_1\cdots b_s}))
=
q^{\#\{i\:|\:a_i < b\}+\#\{j\:|\:b_j <b\}}
L(\substack{a\\b})
M(\substack{a_1\cdots a_r\\b_1\cdots b_s})$.
Now by weight considerations we get that
\begin{align*}
x := 
\omega(L(\substack{a_1\cdots a_r \\ b_1 \cdots b_s}))
&=q^{\#\{i\:|\:a_i > a\}+\#\{j\:|\:b_j >a\}}
L(\substack{a_1\cdots a_r \\ b_1 \cdots b_s}) L(\substack{a\\b})\\
&=
q^{\#\{i\:|\:a_i < b\}+\#\{j\:|\:b_j <b\}}
L(\substack{a\\b})L(\substack{a_1\cdots a_r \\ b_1 \cdots b_s}).
\end{align*}
Using Theorem~\ref{qbar}, we deduce that
\begin{align*}
\overline{x} &=
q^{-\#\{i\:|\:a_i > a\}-\#\{j\:|\:b_j >a\}+r+s
-(\eps_a+\eps_b,\eps_{a_1}+\cdots+\eps_{a_r}+\eps_{b_1}+\cdots+\eps_{b_s})}
L(\substack{a\\b}) 
L(\substack{a_1\cdots a_r\\b_1\cdots b_s})\\
&=q^{\#\{i\:|\:a_i <b\}+\#\{j\:|\:b_j <b\}}
L(\substack{a\\b})L(\substack{a_1\cdots a_r\\b_1\cdots b_s}) = x.
\end{align*}
Hence, $x$ is bar invariant, 
and since it equals $M(\substack{a_1\cdots a_r a \\ b_1\cdots b_s b})$
plus a $q^{-1}\Z[q^{-1}]$-linear combination of other monomials,
we have proved that $x = L(\substack{a_1\cdots a_r a \\ b_1\cdots b_s b})$.
\end{proof}

\begin{thm}\label{pod}
Let $r,s \geq 0$ and $t = \min(r,s)$.
Suppose
$1 \leq a_1,\dots, a_r,b_1,\dots,b_s \leq n$ satisfy the following
property for all $i=1,\dots,t$:
\begin{quote}
If the set $\{a_j - b_k\:|\:i \leq j \leq r, i \leq k \leq s
\text{ such that }a_j > b_k\}$ is non-empty, then $(a_i-b_i)$
is its smallest element.
\end{quote}
Then, 
up to multiplication by a power of $q$,
the dual canonical basis element
$L(\substack{a_1\cdots a_r \\ b_1 \cdots b_s})$ is equal to 
$$
\prod_{\substack{1 \leq i \leq t \\ a_i > b_i}}
(x_{2,a_i} x_{1,b_i} - q^{-1} x_{2,b_i} x_{1,a_i})
\prod_{\substack{1 \leq i \leq t \\ a_i \leq b_i}} x_{2,a_i}x_{1,b_i}
\prod_{t < j \leq r} x_{2,a_j}
\prod_{t < k \leq s} x_{1,b_k}
$$
where the product is taken
in any order.
Every element of the dual canonical basis 
of $\O_q(\M_{2,n})$ can be obtained in this way.
\end{thm}

\begin{proof}
Apply Lemma~\ref{ms} and induction on $t$.
The induction starts from the observation that if
$a_i \leq b_j$ for {\em all} $i,j$ then we have simply that
$L(\substack{a_1\cdots a_r \\ b_1\cdots b_s})
=
M(\substack{a_1\cdots a_r \\ b_1\cdots b_s})$.
\end{proof}

\begin{remark}\rm
Applying $\tau$ to Theorem~\ref{pod}, 
one also obtains a closed formula for the dual canonical basis
of $\O_q(\M_{n,2})$, hence of 
the $U_q(\mathfrak{gl}_2)$-modules $S^\mu(\V_2)$ for all 
$\mu$.
As a special case, we recover the computation by Frenkel and Khovanov
of the dual canonical basis
of the $U_q(\mathfrak{gl}_2)$-module 
$T^d(\V_2)$; see \cite[3.1]{FKK}.
\end{remark}

\begin{remark}\rm
There is one other situation where it is possible to compute 
the canonical/dual canonical bases from $\S$\ref{spowers} explicitly. 
In his thesis, Khovanov also
computed the {\em canonical basis}
of the $U_q(\mathfrak{gl}_2)$-module
$T^d(\V_2)$, which is closely related to the parabolic Kazhdan-Lusztig 
polynomials studied by Lascoux and Sch\"utzenberger in \cite{LS}; 
see \cite[3.4]{FKK}.
The dual statement to this has been derived recently by Cheng, Wang and Zhang
\cite[6.17]{CWZ};
in particular, they give a closed formula for the canonical basis of
$\bw(\V_n) \otimes \bw(\V_n)$, i.e. the canonical basis elements
$K_A^*$ for column strict tableaux $A$ with just {\em two columns}.
\end{remark}

Finally in this section, we want to make precise the 
relationship between the dual canonical basis
of $\O_q(\M_{m,n})$ described here and the dual canonical basis of the 
quantized coordinate algebra $\O_q(\T_{m+n})$
of the group of all upper unitriangular $(m+n)\times (m+n)$-matrices\footnote{Since completing this article, I have learnt of a preprint of
Jakobsen and Zhang \cite{ZJ} which also makes this identification by
similar arguments.}.
Following \cite{BZ}, this is the $\Q(q)$-algebra
on generators $\{t_{i,j}\:|\:1 \leq i < j \leq m+n\}$ subject to the relations
$$
t_{i,k} = \frac{t_{i,j} t_{j,k} - q^{-1} t_{j,k} t_{i,j}}{q-q^{-1}}
$$
for $1 \leq i < j  < k\leq m+n$ and
\begin{align*}
t_{i,j} t_{k,l} &= t_{k,l} t_{i,j}&(i < k, j > l) \text{ or }(i > l)\\
t_{i,j} t_{k,l} &= t_{k,l} t_{i,j} + (q-q^{-1}) t_{i,l} t_{k,j}&(i < k < j < l)\\
t_{i,j} t_{k,j} &= q t_{k,j} t_{i,j}&(i < k)\\
t_{i,j} t_{i,l} &= q t_{i,l} t_{i,j}&(j < l)
\end{align*}
for $1 \leq i < j \leq m+n$ and $1 \leq k < l \leq m+n$.
We view $\O_q(\T_{m+n})$ as an $X_{m+n}$-graded algebra, 
by declaring that the generator $t_{i,j}$ 
is of weight $(\eps_i - \eps_j)$.
In fact, by an observation of Drinfeld proved in \cite{BZ}, 
the algebra $\O_q(\T_{m+n})$ can be identified
with the positive part $\U_{m+n}^+$ of
the quantized enveloping algebra $U_q(\mathfrak{gl}_{m+n})$, so
that $t_{i,i+1}$ is identified with $E_i$ for each $i=1,\dots,m+n-1$.
Under this identification, the bar involution on $\U_{m+n}^+$
also defines a bar involution on 
$\O_q(\T_{m+n})$. Define a different antilinear involution
$\sim$ of $\O_q(\T_{m+n})$ by setting
$\widetilde{x} = 
q^{\frac{1}{2}(\mu,\mu)-\deg(\mu)} \sigma(\overline{x})$
for each $x$ of weight $\mu$.
Here, $\sigma:\O_q(\T_{m+n}) \rightarrow \O_q(\T_{m+n})$ 
is the unique algebra antiautomorphism that fixes the generators
$t_{i,i+1}$ for each $1\leq i < m+n$, and
for a weight $0 \leq \mu \in X_{m+n}$ its degree 
$\deg(\mu)$ is defined from 
$\deg(\eps_i - \eps_{i+1}) = 1$ and 
$\deg(\mu+\nu) = \deg(\mu)+\deg(\nu)$.

To define the dual canonical basis of $\O_q(\T_{m+n})$ following
\cite[$\S$3.5]{LNT}, we must first introduce a PBW basis.
Let $J_{m+n}^d$ denote the set of all
terminal double indexes 
$(\alpha,\beta) \in I_{m+n}^d \times I_{m+n}^d$, such that
$\alpha_i < \beta_i$ for all $i=1,\dots,d$.
For $(\alpha,\beta) \in J_{m+n}^d$,
define 
$$
E^*_{\alpha,\beta} = 
q^{\sum_{i=1}^{m+n}\!\!\nu_i(\nu_i-1)/2}
\,
t_{\alpha_1,\beta_1}\cdots t_{\alpha_d,\beta_d}
$$
where $\nu = \theta(\beta) \in \Lambda_{m+n}$.
This is {\em exactly} the PBW basis element denoted
$\Phi(E^*({\bf m}))$ in \cite{LNT}, parametrized by the multi-segment
${\bf m} = \sum_{i=1}^d [\alpha_i,\beta_i-1]$.
The elements $\{E^*_{\alpha,\beta}\:|\:(\alpha,\beta) \in 
\bigcup_{d \geq 0} J_{m+n}^d\}$ 
give a basis for $\O_q(\T_{m+n})$.
By \cite[3.16]{LNT}, there is for $(\alpha,\beta) \in J_{m+n}^d$ a unique element $G^*_{\alpha,\beta} \in \O_q(\T_{m+n})$ such that
$\widetilde{G^*_{\alpha,\beta}} = G^*_{\alpha,\beta}$ and
\begin{equation}\label{l3}
G^*_{\alpha,\beta} \in E^*_{\alpha,\beta} + \sum_{(\alpha',\beta') \in J_{m+n}^d}
q\Z[q] E^*_{\alpha',\beta'}.
\end{equation}
Moreover, the {\em dual canonical basis} of $\O_q(\T_{m+n})$ is
$\{G^*_{\alpha,\beta}\:|\:(\alpha,\beta) \in \bigcup_{d \geq 0}J_{m+n}^d\}$; it is 
the basis dual to the canonical basis of $\U_{m+n}^+$ under 
a natural bilinear form normalized as in \cite[$\S$3.4]{LNT}.

\begin{thm}
There is an algebra monomorphism
$\phi:\O_q(\M_{m,n}) \rightarrow \O_q(\T_{m+n})$ such that
$\phi(x_{i,j}) = t_{i,j+m}$ for all $1 \leq i \leq m, 1 \leq j \leq n$.
Moreover, given
$\mu \in \Lambda_m$, $\nu \in \Lambda_n$ with $|\mu|=|\nu| = d$
and any $(\alpha,\beta) \in (I_\mu\times I_\nu)^-$, we have that
$$
\phi(M_{\alpha,\beta}) = q^{-\sum_{i=1}^n \nu_i(\nu_i-1)/2} E^*_{\alpha,\beta'},
\qquad
\phi(L_{\alpha,\beta}) = q^{-\sum_{i=1}^n \nu_i(\nu_i-1)/2} G^*_{\alpha,\beta'},
$$
where $\beta' = (\beta_1+m,\dots,\beta_d+m)$.
\end{thm}

\begin{proof}
It is clear from the relations that $\phi$ is a well-defined algebra homomorphism. Also, it sends $M_{\alpha,\beta}$ to 
$q^{-\sum_{i=1}^n \nu_i(\nu_i-1)/2} E^*_{\alpha,\beta'}$, hence it is injective.
It just remains to show that it sends $L_{\alpha,\beta}$ to
$q^{-\sum_{i=1}^n \nu_i(\nu_i-1)/2} G^*_{\alpha,\beta'}$. This follows easily
comparing (\ref{l2}) and (\ref{l3}) as soon as we have checked that
$\phi(q^{\sum_{i=1}^n \nu_i(\nu_i-1)/2}L_{\alpha,\beta})$ is invariant
under the antilinear involution $\sim$.
One checks from the definition that
$\widetilde{t_{i,j}} = t_{i,j}$ for all $1 \leq i < j \leq m+n$, and that
$\widetilde{xy} = q^{(\mu,\nu)} \widetilde{y}\, \widetilde{x}$
for all $x,y$ of weights $\mu,\nu$, respectively.
Combining this with Theorem~\ref{qbar}, it follows by induction on degree that
$$
\phi(\overline{x}) = q^{-\sum_{i=1}^n \nu_i(\nu_i-1)} \widetilde{\phi(x)}
$$
for any $x \in \O_q(\M_{m,n})_{\mu,\nu}$.
Hence, 
$\widetilde{\phi(L_{\alpha,\beta})} = 
q^{\sum_{i=1}^n \nu_i(\nu_i-1)} \phi(L_{\alpha,\beta})$.
\end{proof}

\begin{remark}\label{positivity}\rm
This theorem means that one can appeal to the extensive literature 
on dual canonical bases of $\O_q(\T_{m+n})$ in order to obtain
powerful results about the dual canonical basis of 
$\O_q(\M_{m,n})$ too. For example, by dualizing 
\cite[14.4.13(b)]{Lubook}, it follows that the structure constants 
for multiplication in $\O_q(\M_{m,n})$ relative to the dual canonical basis
in fact all lie in $\N[q,q^{-1}]$.
\end{remark}

\section{Polynomial representations}\label{spoly}

Assume throughout the section that $m \leq n$ and
that $\mu \in \Lambda_l$ is a weight with
$|\mu| = d$, such that the conjugate partition
$\lambda = \mu'$ lies in $\Lambda_m^+$.
Recall the definitions from $\S$\ref{spowers}
of the $\U_n$-modules
$\bw^{\mu}(\V_n)$ and $S^{\lambda}(\V_n)$. 
The $\nu$-weight space of the first one
has the two natural bases
$N_A$ and $K_A$ parametrized by $\Col(\mu,\nu)$,
while the $\nu$-weight space of the second one has the two
natural bases $M_B$ and $L_B$ parametrized by $\Row(\lambda,\nu)$.
By the Littlewood-Richardson rule, the $\U_n$-module
$\bw^\mu(\V_n)$ resp. $S^\lambda(\V_n)$ 
has a composition factor of highest weight $\lambda$
appearing with multiplicity one, and all the other composition factors
are of highest weight $< \lambda$ resp. $> \lambda$ in the dominance ordering.
Hence, the space $\hom_{\U_n}(\bw^\mu(\V_n), S^\lambda(\V_n))$ 
is one dimensional. 
We define $P^\lambda(\V_n)$ to be the image
of any non-zero 
homomorphism $\bw^\mu(\V_n) \rightarrow S^\lambda(\V_n)$.
This is the well known realization of the irreducible
polynomial representation of
$\U_n$ of highest weight $\lambda$ as a submodule
of $S^\lambda(\V_n)$.
We should note that since the spaces $\bw^\mu(\V_n)$ are 
at least {\em isomorphic} 
for all $\mu$ with $\mu' = \lambda$,
the one dimensionality of $\hom_{\U_n}(\bw^\mu(\V_n), S^\lambda(\V_n))$
implies that $P^\lambda(\V_n)$ is always the {\em same} subspace of
$S^\lambda(\V_n)$, independent of the particular choice of $\mu$.

Let us write down a canonical
generator for the space
$\hom_{\U_n}(\bw^\mu(\V_n), S^\lambda(\V_n))$.
To do this, we identify $S^\lambda(\V_n)$
with a one-sided weight space of $\O_q(\M_{m,n})$ according to 
Theorem~\ref{qca}.
Given $\beta \in I_n^d$ with $\beta_1 < \cdots < \beta_d$, 
define the {\em quantum flag minor}
\begin{align}\label{qfm}
D_\beta &= \sum_{w \in S_d} (-q)^{\ell(w)} x_{w1, \beta_1} x_{w2,\beta_2}
\cdots x_{wd, \beta_d}\\ &=
\sum_{w \in S_d} (-q)^{-\ell(w)} x_{d, \beta_{wd}} \cdots x_{2,\beta_{w2}}
x_{1, \beta_{w1}}.\label{qfm2}
\end{align}
Recalling Theorem~\ref{qbar},
it is immediate from this definition that
$D_\beta$ is bar invariant,
hence it coincides with the dual canonical basis element
$L_{\alpha,\beta}$ where $\alpha = (1,2,\dots,d)$. Now for $A \in \Col(\mu,\nu)$,
define
\begin{equation}\label{vad}
V_A := D_{\alpha_1} D_{\alpha_2} \cdots D_{\alpha_d},
\end{equation}
where $\alpha_i$ denotes the multi-index obtained by reading the entries in the $i$th column of $A$ from bottom to top. Thus, $V_A$ is the product of the quantum flag minors corresponding to the columns of the tableau $A$.
Clearly it belongs to the one-sided weight space
$S^\lambda(\V_n)$ 
of $\O_q(\M_{m,n})$, so we can define a linear map
\begin{equation}
\xi_\mu:\bw^\mu(\V_n) \rightarrow S^\lambda(\V_n)
\end{equation}
by setting
$\xi_\mu(N_A) = V_A$
for all $A \in \Col(\mu,\nu)$.
We note finally that if $A$ is the unique element of $\Col(\mu,\lambda)$,
so all entries on the $i$th row of $A$ are equal to $i$, then
\begin{equation}\label{leading}
V_A = M_{R(A)} + \text{(a $\Z[q,q^{-1}]$-linear combination of $M_B$'s for
$B < R(A)$).}
\end{equation}
Of course, the rectification 
$R(A)$ in this case is just the tableau of row shape $\lambda$ having all
entries on its $i$th row equal to $i$.
The proof of (\ref{leading}) is a straightforward 
consequence of the defining relations in $\O_q(\M_{m,n})$.

\begin{lem}\label{quicker}
The map $\xi_\mu$ is a non-zero $\U_n$-module 
homomorphism. 
\iffalse
Moreover, 
$\xi_\mu(\overline{v}) = \overline{\xi_\mu(v)}$ for all
$v \in \bw^\mu(\V_n)$.
\fi
\end{lem}

\begin{proof}
For each $i$, identify $T^{\mu_i}(\V_n)$ with
a submodule of $\O_q(\M_{m,n})$ by identifying
$v_{\alpha_1}\otimes\cdots\otimes v_{\alpha_{\mu_i}}$
with $x_{\mu_i,\alpha_1} x_{\mu_i-1,\alpha_2}\cdots x_{1,\alpha_{\mu_i}}$.
In this way, 
$T^d(\V_n) = T^{\mu_1}(\V_n) \otimes\cdots\otimes T^{\mu_l}(\V_n)$
is identified with a submodule of $\O_q(\M_{m,n})^{\otimes l}$.
Let $A \in \Col(\mu,\nu)$ be a column strict tableau, and let
$\alpha_i$ denote the multi-index obtained by reading the entries in the $i$th
column of $A$ from bottom to top.
Comparing (\ref{ya}) with the right hand side of
(\ref{qfm2}), 
the basis element $N_A = Y_{\alpha_1}\otimes\cdots\otimes 
Y_{\alpha_l}$ of $\bw^\mu(\V_n) \subseteq T^d(\V_n)$ 
corresponds under this identification
to the tensor product of quantum flag minors
$D_{\alpha_1} \otimes\cdots\otimes D_{\alpha_l} \in \O_q(\M_{m,n})^{\otimes l}$.
Since $\O_q(\M_{m,n})$ is a polynomial $\U_n$-algebra, 
multiplication defines a
$\U_n$-module homomorphism
$\O_q(\M_{m,n})^{\otimes l} \rightarrow \O_q(\M_{m,n})$ mapping
$N_A$ to $V_A$.
Hence, $\xi_\mu$ is a $\U_n$-module homomorphism, and it
is non-zero by (\ref{leading}).
\iffalse
Finally, 
to see that it commutes with the bar involution, note that
$- \circ \xi_\mu \circ -$ is another 
$\U_n$-module homomorphism from $\bw^\mu(\V_n)$ to $S^\lambda(\V_n)$,
hence it is equal to $\xi_\mu$ up to a scalar by the one dimensionality
of the space $\hom_{\U_n}(\bw^\mu(\V_n), S^\lambda(\V_n))$.
To show that the scalar is $1$, note that $\overline{N_A} = N_A$,
while by (\ref{leading}) the $M_{R(A)}$-coefficient of both $V_A$ and $\overline{V_A}$ equals $1$.
\fi
\end{proof}

\begin{thm}\label{main1}
For any $\nu \in \Lambda_n$ and $A \in \Col(\mu,\nu)$, we have that
$$
\xi_{\mu}(K_A) = \left\{
\begin{array}{ll}
L_{R(A)}&\text{if $A \in \Std(\mu,\nu)$,}\\
0&\text{otherwise.}
\end{array}\right.
$$
The vectors
$\{V_A\:|\:A \in \Std(\mu,\nu)\}$ and
$\{L_B\:|\:B \in \Dom(\lambda,\nu)\}$
give natural bases for the $\nu$-weight space $P_\nu^\lambda(\V_n)$
of $P^\lambda(\V_n)$. Moreover, 
for $A \in \Col(\mu,\nu)$,
we have that
$$
V_A 
= \sum_{B \in \Std(\mu,\nu)} k_{A,B}^*(q^{-1}) L_{R(B)}.
$$
\end{thm}

\begin{proof}
Recall the subring $\A_\infty$ of $\Q(q)$ from $\S$\ref{skl}.
Let $\bw^\mu(\V_n)_\infty$ resp. $S^\lambda(\V_n)_\infty$ be the
$\A_\infty$-submodule of $\bw^\mu(\V_n)$ resp. $S^\lambda(\V_n)$
generated by all the $N_A$'s resp. $M_A$'s.
It is an upper crystal lattice at $q=\infty$ in the sense of \cite{KaG}, and
the images of the $N_A$'s resp. $M_A$'s
in
$\bw^\mu(\V_n)_\infty / q^{-1} \bw^\mu(\V_n)_\infty$
resp.
$S^\lambda(\V_n)_\infty / q^{-1} S^\lambda(\V_n)_\infty$ form an upper crystal
base at $q=\infty$.
The action of the upper crystal operators on this upper crystal base
is described by the crystal
$\bigcup \Col(\mu,\nu)$ resp.
$\bigcup \Row(\lambda,\nu)$
from $\S$\ref{scombinatorics}.
Finally, $(\Q \otimes_{\Z} \bw^\mu(\VV_n), \overline{\bw^\mu(\V_n)_\infty},
\bw^\mu(\V_n)_\infty)$ resp.
$(\Q \otimes_{\Z} S^\lambda(\VV_n), \overline{S^\lambda(\V_n)_\infty},
S^\lambda(\V_n)_\infty)$ is a balanced triple, and the dual canonical basis
of $\bw^\mu(\V_n)$ resp. $S^\lambda(\V_n)$
is the corresponding lift of the upper crystal base.
This puts us in the setup of \cite[$\S$5]{KaG}.

Take any $\nu \in \Lambda_n$ and $A \in \Col(\mu,\nu)$ such that
$\tilde e_i (A) = \varnothing$ for all $i$.
Then, by \cite[5.1.1]{KaG}, $K_A$ is a non-zero highest weight vector
in $\bw^\mu(\V_n)$ of weight $\nu$. 
Since all composition factors of
$\bw^\mu(\V_n)$ are of highest weight $\leq \lambda$,
we have that
$\nu \leq \lambda$.
Since all composition factors of $S^\lambda(\V_n)$ are of highest weight
$\geq \lambda$, we deduce that
$\xi_\mu(K_A) = 0$ unless in fact $\nu = \lambda$. 
In that case,
there is only one tableau $A \in \Col(\mu,\lambda)$,
and so we must have that $K_A = N_A$.
Since $\tilde e_i (R(A)) = \varnothing$ for all $i$ too, we get 
by \cite[5.1.1]{KaG} once more that
$L_{R(A)}$ is a highest weight vector in $S^\lambda(\V_n)$ of weight
$\lambda$. Hence
$\xi_\mu(K_A) = V_A = c L_{R(A)}$ for some non-zero scalar $c \in \Q(q)$.
Since $L_{R(A)} = M_{R(A)} + $(a $q^{-1}\Z[q^{-1}]$-linear combination of
$M_B$'s for $B < R(A)$) we deduce from (\ref{leading}) that $c=1$.
Hence, $\xi_\mu(K_A) = L_{R(A)}$ in this special case.

Now for the general case, the point is that
there are two possibly different balanced triples
in $P^\lambda(\V_n)$, one arising as a quotient of
the balanced triple $(\Q \otimes_{\Z} \bw^\mu(\VV_n), \overline{\bw^\mu(\V_n)_\infty},
\bw^\mu(\V_n)_\infty)$, the other arising from the intersection
with the balanced triple
$(\Q \otimes_{\Z} S^\lambda(\VV_n), \overline{S^\lambda(\V_n)_\infty},
S^\lambda(\V_n)_\infty)$.
We have just checked in the previous paragraph that these two balanced
triples agree on the highest weight space of the irreducible
module $P^\lambda(\V_n)$. Hence by \cite[5.2.2]{KaG}, they agree everywhere.
This shows in particular that the map $\xi_\mu$ maps the upper crystal lattice
$\bw^\mu(\V_n)_\infty$ into $S^\lambda(\V_n)_\infty$, so we get an induced map 
$\bar\xi_\mu:\bw^\mu(\V_n)_\infty / q^{-1} \bw^\mu(\V_n)_\infty
\rightarrow S^\lambda(\V_n)_\infty /q^{-1} S^\lambda(\V_n)_\infty$
commuting with the actions of the upper crystal operators.
Moreover, the following diagram commutes
$$
\begin{CD}
\Q\otimes_\Z\bw^\mu(\VV_n)&@>\sim>>&\bw^\mu(\V_n)_\infty / q^{-1} \bw^\mu(\V_n)_\infty
\\
@V\xi_\mu VV&&@VV\bar\xi_\mu V\\
\Q\otimes_\Z S^\lambda(\VV_n)&@>\sim>>&S^\lambda(\V_n)_\infty / q^{-1}S^\lambda(\V_n)_\infty
\end{CD}
$$
where the top and bottom maps are the canonical isomorphisms arising from the
balanced triples.
It now suffices to complete the proof of the first statement of the theorem
to verify it at the level of local crystal bases.
If $A \in \Col(\mu,\nu)$ 
satisfies $\tilde e_i(A) = \varnothing$ for all $i$, we are done by the
previous paragraph. The general case follows by applying crystal operators,
recalling the characterization of the set $\bigcup\Std(\mu,\nu)$
and the map $R$ in terms of crystals from $\S$\ref{scombinatorics}.

It follows immediately that 
$\{L_A\:|\:A \in \Dom(\lambda,\nu)\}$ is a basis for
$P^\lambda_\nu(\V_n)$.
By (\ref{ekl4}) and Corollary~\ref{kinv}, we have 
for any $A \in \Col(\mu,\nu)$ that
$$
N_A = 
\sum_{B \in \Col(\mu,\nu)}
k_{A,B}^*(q^{-1})  K_B.
$$
Applying the map $\xi_{\mu}$, we get the formula for $V_A$.
Finally unitriangularity of the transition matrix
implies that $\{V_A\:|\:A \in \Std(\mu,\nu)\}$ is also a basis
for $P^\lambda_\nu(\V_n)$.
\end{proof}

\begin{remark}\label{tm}\rm
The basis $\{L_A\:|\:A \in \bigcup_{\nu}\Dom(\lambda,\nu)\}$ 
for $P^\lambda(\V_n)$ 
is Kashiwara's upper global crystal base (by the proof of Theorem~\ref{main1})
or Lusztig's {dual canonical basis} (by
Remarks~\ref{eq1} and \ref{eq2} below).
It is the same basis 
independent of the choice of $\mu$. On the other hand, the 
basis $\{V_A\:|\:A \in \bigcup_{\nu}
\Std(\mu,\nu)\}$ definitely does depend on 
$\mu$. Thus, we obtain a family of standard monomial bases for $P^\lambda(\V_n)$,
one for each $\mu$ with $\mu' = \lambda$. These bases are not new; 
for instance,
they were already constructed
in \cite[4.4]{LT} by a similar approach to the one here.
In the case that $\mu$ is itself a partition,
this basis is the $q$-analogue of the
classical standard monomial basis.
Note finally by the definition (\ref{vad}) and Remark~\ref{positivity} that
the coefficients of the polynomials $k_{A,B}^*(q)$ appearing
in Theorem~\ref{main1} are non-negative integers.
\end{remark}

\begin{example}\rm
We list the polynomials $k_{A,B}^*(q)$
for $\mu = (3,2,2,1)$,
$\nu = (2,2,2,1,1)$
and all $A,B \in \Std(\mu,\nu)$, i.e. part of the transition matrix
from the standard monomial to the dual canonical basis of $P^\lambda(\V_n)$,
where $\lambda = (4,3,1)$.
We pick this example in order to point out that
the $AB$-entry of this matrix is the same as the $AB$-entry 
of the matrix computed by Leclerc and Toffin in \cite{LTo}; in particular 
\cite{LTo}
gives a simple algorithm to compute these polynomials.
$$
\begin{array}{l|ccccccccccccc}
\!\!\!\!\!\scriptstyle k^*_{A,B}(q)&\!\!\!\!\begin{array}{llll}
\scriptscriptstyle 3\vspace{-2.6mm}\\
\scriptscriptstyle 224\vspace{-2.6mm}\\
\scriptscriptstyle 1135\vspace{-0.1mm}
\end{array}
&\!\!\!\!\!\begin{array}{llll}
\scriptscriptstyle 3\vspace{-2.6mm}\\
\scriptscriptstyle 225\vspace{-2.6mm}\\
\scriptscriptstyle 1134\vspace{-0.1mm}
\end{array}
&\!\!\!\!\!\begin{array}{llll}
\scriptscriptstyle 3\vspace{-2.6mm}\\
\scriptscriptstyle 234\vspace{-2.6mm}\\
\scriptscriptstyle 1125\vspace{-0.1mm}
\end{array}
&\!\!\!\!\!\begin{array}{llll}
\scriptscriptstyle 3\vspace{-2.6mm}\\
\scriptscriptstyle 235\vspace{-2.6mm}\\
\scriptscriptstyle 1124\vspace{-0.1mm}
\end{array}
&\!\!\!\!\!\begin{array}{llll}
\scriptscriptstyle 3\vspace{-2.6mm}\\
\scriptscriptstyle 245\vspace{-2.6mm}\\
\scriptscriptstyle 1123\vspace{-0.1mm}
\end{array}
&\!\!\!\!\!\begin{array}{llll}
\scriptscriptstyle 4\vspace{-2.6mm}\\
\scriptscriptstyle 225\vspace{-2.6mm}\\
\scriptscriptstyle 1133\vspace{-0.1mm}
\end{array}
&\!\!\!\!\!\begin{array}{llll}
\scriptscriptstyle 4\vspace{-2.6mm}\\
\scriptscriptstyle 233\vspace{-2.6mm}\\
\scriptscriptstyle 1125\vspace{-0.1mm}
\end{array}
&\!\!\!\!\!\begin{array}{llll}
\scriptscriptstyle 4\vspace{-2.6mm}\\
\scriptscriptstyle 235\vspace{-2.6mm}\\
\scriptscriptstyle 1123\vspace{-0.1mm}
\end{array}
&\!\!\!\!\!\begin{array}{llll}
\scriptscriptstyle 5\vspace{-2.6mm}\\
\scriptscriptstyle 224\vspace{-2.6mm}\\
\scriptscriptstyle 1133\vspace{-0.1mm}
\end{array}
&\!\!\!\!\!\begin{array}{llll}
\scriptscriptstyle 5\vspace{-2.6mm}\\
\scriptscriptstyle 233\vspace{-2.6mm}\\
\scriptscriptstyle 1124\vspace{-0.1mm}
\end{array}
&\!\!\!\!\!\begin{array}{llll}
\scriptscriptstyle 5\vspace{-2.6mm}\\
\scriptscriptstyle 234\vspace{-2.6mm}\\
\scriptscriptstyle 1123\vspace{-0.1mm}
\end{array}
&\!\!\!\!\!\begin{array}{llll}
\scriptscriptstyle 4\vspace{-2.6mm}\\
\scriptscriptstyle 335\vspace{-2.6mm}\\
\scriptscriptstyle 1122\vspace{-0.1mm}
\end{array}
&\!\!\!\!\!\begin{array}{llll}
\scriptscriptstyle 5\vspace{-2.6mm}\\
\scriptscriptstyle 334\vspace{-2.6mm}\\
\scriptscriptstyle 1122\vspace{-0.1mm}
\end{array}\\\hline
\begin{array}{llll}
\scriptscriptstyle 3\vspace{-2.6mm}\\
\scriptscriptstyle 224\vspace{-2.6mm}\\
\scriptscriptstyle 1135\vspace{-1.2mm}
\end{array}\!\!&\scriptstyle 1&\scriptstyle \cdot&\scriptstyle \cdot&\scriptstyle \cdot&\scriptstyle \cdot&\scriptstyle \cdot&\scriptstyle \cdot&\scriptstyle \cdot&\scriptstyle \cdot&\scriptstyle \cdot&\scriptstyle \cdot&\scriptstyle \cdot&\scriptstyle \cdot\\
\begin{array}{llll}
\scriptscriptstyle 3\vspace{-2.6mm}\\
\scriptscriptstyle 225\vspace{-2.6mm}\\
\scriptscriptstyle 1134\vspace{-1.2mm}
\end{array}\!\!&\scriptstyle q&\scriptstyle 1&\scriptstyle \cdot&\scriptstyle \cdot&\scriptstyle \cdot&\scriptstyle \cdot&\scriptstyle \cdot&\scriptstyle \cdot&\scriptstyle \cdot&\scriptstyle \cdot&\scriptstyle \cdot&\scriptstyle \cdot&\scriptstyle \cdot\\
\begin{array}{llll}
\scriptscriptstyle 3\vspace{-2.6mm}\\
\scriptscriptstyle 234\vspace{-2.6mm}\\
\scriptscriptstyle 1125\vspace{-1.2mm}
\end{array}\!\!&\scriptstyle q&\scriptstyle \cdot&\scriptstyle 1&\scriptstyle \cdot&\scriptstyle \cdot&\scriptstyle \cdot&\scriptstyle \cdot&\scriptstyle \cdot&\scriptstyle \cdot&\scriptstyle \cdot&\scriptstyle \cdot&\scriptstyle \cdot&\scriptstyle \cdot\\
\begin{array}{llll}
\scriptscriptstyle 3\vspace{-2.6mm}\\
\scriptscriptstyle 235\vspace{-2.6mm}\\
\scriptscriptstyle 1124\vspace{-1.2mm}
\end{array}\!\!&\scriptstyle q^2&\scriptstyle q&\scriptstyle q&\scriptstyle 1&\scriptstyle \cdot&\scriptstyle \cdot&\scriptstyle \cdot&\scriptstyle \cdot&\scriptstyle \cdot&\scriptstyle \cdot&\scriptstyle \cdot&\scriptstyle \cdot&\scriptstyle \cdot\\
\begin{array}{llll}
\scriptscriptstyle 3\vspace{-2.6mm}\\
\scriptscriptstyle 245\vspace{-2.6mm}\\
\scriptscriptstyle 1123\vspace{-1.2mm}
\end{array}\!\!&\scriptstyle q&\scriptstyle q^2&\scriptstyle q^2&\scriptstyle q&\scriptstyle 1&\scriptstyle \cdot&\scriptstyle \cdot&\scriptstyle \cdot&\scriptstyle \cdot&\scriptstyle \cdot&\scriptstyle \cdot&\scriptstyle \cdot&\scriptstyle \cdot\\
\begin{array}{llll}
\scriptscriptstyle 4\vspace{-2.6mm}\\
\scriptscriptstyle 225\vspace{-2.6mm}\\
\scriptscriptstyle 1133\vspace{-1.2mm}
\end{array}\!\!&\scriptstyle q&\scriptstyle q^2&\scriptstyle \cdot&\scriptstyle \cdot&\scriptstyle \cdot&\scriptstyle 1&\scriptstyle \cdot&\scriptstyle \cdot&\scriptstyle \cdot&\scriptstyle \cdot&\scriptstyle \cdot&\scriptstyle \cdot&\scriptstyle \cdot\\
\begin{array}{llll}
\scriptscriptstyle 4\vspace{-2.6mm}\\
\scriptscriptstyle 233\vspace{-2.6mm}\\
\scriptscriptstyle 1125\vspace{-1.2mm}
\end{array}\!\!&\scriptstyle q&\scriptstyle \cdot&\scriptstyle q^2&\scriptstyle \cdot&\scriptstyle \cdot&\scriptstyle \cdot&\scriptstyle 1&\scriptstyle \cdot&\scriptstyle \cdot&\scriptstyle \cdot&\scriptstyle \cdot&\scriptstyle \cdot&\scriptstyle \cdot\\
\begin{array}{llll}
\scriptscriptstyle 4\vspace{-2.6mm}\\
\scriptscriptstyle 235\vspace{-2.6mm}\\
\scriptscriptstyle 1123\vspace{-1.2mm}
\end{array}\!\!&\scriptstyle 2q^2&\scriptstyle q^3&\scriptstyle q^3&\scriptstyle q^2&\scriptstyle q&\scriptstyle q&\scriptstyle q&\scriptstyle 1&\scriptstyle \cdot&\scriptstyle \cdot&\scriptstyle \cdot&\scriptstyle \cdot&\scriptstyle \cdot\\
\begin{array}{llll}
\scriptscriptstyle 5\vspace{-2.6mm}\\
\scriptscriptstyle 224\vspace{-2.6mm}\\
\scriptscriptstyle 1133\vspace{-1.2mm}
\end{array}\!\!&\scriptstyle q^2&\scriptstyle \cdot&\scriptstyle \cdot&\scriptstyle \cdot&\scriptstyle \cdot&\scriptstyle q&\scriptstyle \cdot&\scriptstyle \cdot&\scriptstyle 1&\scriptstyle \cdot&\scriptstyle \cdot&\scriptstyle \cdot&\scriptstyle \cdot\\
\begin{array}{llll}
\scriptscriptstyle 5\vspace{-2.6mm}\\
\scriptscriptstyle 233\vspace{-2.6mm}\\
\scriptscriptstyle 1124\vspace{-1.2mm}
\end{array}\!\!&\scriptstyle q^2&\scriptstyle q&\scriptstyle q^3&\scriptstyle q^2&\scriptstyle \cdot&\scriptstyle \cdot&\scriptstyle q&\scriptstyle \cdot&\scriptstyle \cdot&\scriptstyle 1&\scriptstyle \cdot&\scriptstyle \cdot&\scriptstyle \cdot\\
\begin{array}{llll}
\scriptscriptstyle 5\vspace{-2.6mm}\\
\scriptscriptstyle 234\vspace{-2.6mm}\\
\scriptscriptstyle 1123\vspace{-1.2mm}
\end{array}\!\!&\scriptstyle 2q^3&\scriptstyle q^2&\scriptstyle q^4&\scriptstyle q^3&\scriptstyle q^2&\scriptstyle q^2&\scriptstyle q^2&\scriptstyle q&\scriptstyle q&\scriptstyle q&\scriptstyle 1&\scriptstyle \cdot&\scriptstyle \cdot\\
\begin{array}{llll}
\scriptscriptstyle 4\vspace{-2.6mm}\\
\scriptscriptstyle 335\vspace{-2.6mm}\\
\scriptscriptstyle 1122\vspace{-1.2mm}
\end{array}\!\!&\scriptstyle q^3&\scriptstyle \cdot&\scriptstyle \cdot&\scriptstyle \cdot&\scriptstyle q^2&\scriptstyle q^4&\scriptstyle q^2&\scriptstyle q^3+q&\scriptstyle \cdot&\scriptstyle \cdot&\scriptstyle \cdot&\scriptstyle 1&\scriptstyle \cdot\\
\begin{array}{llll}
\scriptscriptstyle 5\vspace{-2.6mm}\\
\scriptscriptstyle 334\vspace{-2.6mm}\\
\scriptscriptstyle 1122
\end{array}\!\!&\scriptstyle q^4&\scriptstyle q^3&\scriptstyle \cdot&\scriptstyle q^2&\scriptstyle q^3&\scriptstyle q^5&\scriptstyle q^3&\scriptstyle q^4+q^2&\scriptstyle q^4&\scriptstyle q^2&\scriptstyle q^3+q&\scriptstyle q&\scriptstyle 1\\
\end{array}
$$
\end{example}

Now let us prepare to dualize.
Recall the spaces $\widetilde S^\lambda(\V_n)$ and
$\widetilde\bw^\mu(\V_n)$ from $\S$\ref{spowers}, and the non-degenerate
pairings $(.,.)$ between $S^\lambda(\V_n)$ and $\widetilde S^\lambda(\V_n)$
and between $\bw^\mu(\V_n)$ and $\widetilde\bw^\mu(\V_n)$.
The space $\hom_{\U_n}(\widetilde S^\lambda(\V_n), \widetilde \bw^\mu(\V_n))$
is also one dimensional, and a canonical generator is given by the map
\begin{equation}
\xi^*_\mu:\widetilde S^\lambda(\V_n) \rightarrow \widetilde \bw^\mu(\V_n)
\end{equation}
that is dual to $\xi_\mu$ in the sense that
$(v, \xi^*_\mu(w)) = (\xi_\mu(v),w)$ for all
$v \in \bw^\mu(\V_n), w \in \widetilde S^\lambda(\V_n)$.
Define $\widetilde P^\lambda(\V_n)$ to be the cokernel of $\xi^*_\mu$
(or indeed any non-zero homomorphism $\widetilde S^\lambda(\V_n)
\rightarrow \widetilde\bw^\mu(\V_n)$). This is another realization
of the irreducible polynomial representation of $\U_n$ as a quotient of
$\widetilde S^\lambda(\V_n)$. It is always the {\em same}
quotient of $\widetilde S^\lambda(\V_n)$ independent of the particular choice
of $\mu$. Actually in practise we will often 
view $\widetilde P^\lambda(\V_n)$
as a submodule of $\widetilde\bw^\mu(\V_n)$ via the map $\xi^*_\mu$,
though of course this identification does depend on 
our fixed
choice of
$\mu$. The pairing 
$(.,.)$ 
between $S^\lambda(\V_n)$ and $\widetilde S^\lambda(\V_n)$ induces
a well-defined non-degenerate pairing
\begin{equation}\label{ppair}
(.,.):P^\lambda(\V_n) \times \widetilde P^\lambda(\V_n) \rightarrow \Q(q).
\end{equation}
Finally, for any $A \in \Row(\lambda,\nu)$, define
$V_A^* = \xi^*_\mu(M_A^*)
\in \widetilde P^\lambda(\V_n)$.

\begin{thm}\label{main2}
For any $\nu \in \Lambda_n$ and $A \in \Row(\lambda,\nu)$, we have that
$$
\xi^*_{\mu}(L_A^*) =
\left\{
\begin{array}{ll}
K_{R^{-1}(A)}^*&\text{if $A \in \Dom(\lambda,\nu)$,}\\
0&\text{otherwise.}
\end{array}\right.
$$
The vectors
$\{V_A^*\:|\:A \in \Dom(\lambda,\nu)\}$
 and
$\{K_B^*\:|\:B \in \Std(\mu,\nu)\}$
give two natural bases for the $\nu$-weight space
$\widetilde P^\lambda_\nu(\V_n)$ of $\widetilde P^\lambda(\V_n)$.
Moreover, 
$(L_A, K_B^*) = \delta_{A,R(B)}$
for $A \in \Dom(\lambda,\nu)$,
$B \in \Std(\mu,\nu)$.
Finally, for any $A \in \Row(\lambda,\nu)$, we have that
$$
V_A^* = \sum_{B \in \Dom(\mu,\nu)} l_{A,B}(q^{-1}) K_{R^{-1}(B)}^*.
$$
\end{thm}

\begin{proof}
For $A \in \Col(\mu,\nu)$, $B \in \Row(\lambda,\nu)$,
$(K_A, \xi_\mu^*(L_B^*)) = (\xi_\mu(K_A), L_B^*)$, which by
Theorems~\ref{bermuda} and \ref{main1} is zero unless $A \in \Std(\mu,\nu)$ and
$B = R(A)$. 
Now argue as in the last paragraph of the proof of 
Theorem~\ref{main1} to get the remaining
statements.
\end{proof}

\begin{remark}\label{eq1}\rm
Note that (unlike in Theorem~\ref{main1}) 
the particular choice of $\mu$ here is irrelevant: it only affects the 
parametrization of the bases not the bases themselves, so one may as well
take $\mu = \lambda'$.
Using Lusztig's results \cite[27.1.7,27.2.4]{Lubook} on filtrations of based modules, it is
not hard to prove Theorem~\ref{main2} directly,
instead of by dualizing Theorem~\ref{main1}. 
This identifies the basis 
$\{K_A^*\:|\:A \in \bigcup_{\nu}\Std(\mu,\nu)\}$
for $\widetilde P^\lambda(\V_n)$ directly with 
the {\em canonical basis} in the sense
of Lusztig, which is the {\em lower global crystal base}
of Kashiwara (by Remarks~\ref{tm} and \ref{eq2}).
The basis $\{V_A^*\:|\:A \in \bigcup_{\nu}
\Dom(\lambda,\nu)\}$ is the {\em semi-standard basis}
of Dipper and James \cite{DJ2}. 
\end{remark}

\begin{remark}\label{eq2}\rm
We proved in Theorem~\ref{main2} that the basis
$\{L_A\:|\:A \in \bigcup_{\nu} \Dom(\lambda,\nu)\}$
for $P^\lambda(\V_n)$ is dual to the basis
$\{K_A^*\:|\:A \in \bigcup_{\nu}\Std(\mu,\nu)\}$
for $\widetilde P^\lambda(\V_n)$ under the pairing $(.,.)$ from (\ref{ppair}).
We can give a more familiar definition of this pairing as follows.
Let $A \in \Col(\mu,\lambda)$ be the tableau having all entries in
its $i$th row equal to $i$. 
Then, $V_A = L_{R(A)}$
and $V_{R(A)}^* = K_A^* $ are the canonical highest weight vectors in
$P^\lambda(\V_n)$ and $\widetilde P^\lambda(\V_n)$, respectively.
By Theorem~\ref{main2}, we have that $(V_A, V_{R(A)}^*) = 1$.
The pairing $(.,.)$ is characterized uniquely by this property
and the fact from Lemma~\ref{mis} that
$(uv,w) = (v,\tau(u) w)$
for all $u \in \U_n, v \in P^\lambda(\V_n),
w \in \widetilde P^\lambda(\V_n)$.
\end{remark}

\begin{remark}\rm
The constructions in this section actually yield bases 
for the $\Z[q,q^{-1}]$-forms
$P^\lambda(\VV_n)$ and $\widetilde P^\lambda(\VV_n)$, meaning
the image resp. cokernel of the restriction of the map
$\xi_\mu$ resp. $\xi_\mu^*$ to $\bw^\mu(\VV_n)$ resp.
$S^\lambda(\VV_n)$.
It is only here that the essential difference between the two 
constructions shows up:
$\widetilde P^\lambda(\VV_n)$ is the $\Z[q,q^{-1}]$-lattice in $\widetilde P^\lambda(\V_n)$
obtained by applying $\UU_n$ to the canonical highest weight vector from
Remark~\ref{eq2},
and $P^\lambda(\VV_n)$ is the dual lattice under the pairing
$(.,.)$.
\end{remark}

\end{document}